\documentclass{amsart}

\usepackage{amsmath}
\usepackage{amscd}
\usepackage{amssymb}
\newcommand{\cal}{\mathcal}
\newcommand{\bk}{{\bf k}}
\newcommand{\bs}{{\bf s}}
\newcommand{\bC}{{\Bbb C}}
\newcommand{\bE}{{\Bbb E}}
\newcommand{\bL}{{\Bbb L}}
\newcommand{\bN}{{\Bbb N}}
\newcommand{\bP}{{\Bbb P}}
\newcommand{\bQ}{{\Bbb Q}}
\newcommand{\bR}{{\Bbb R}}
\newcommand{\bZ}{{\Bbb Z}}
\newcommand{\cA}{{\cal A}}
\newcommand{\cB}{{\cal B}}
\newcommand{\cC}{{\cal C}}
\newcommand{\cD}{{\cal D}}
\newcommand{\cE}{{\cal E}}
\newcommand{\cF}{{\cal F}}
\newcommand{\cG}{{\cal G}}
\newcommand{\cH}{{\cal H}}
\newcommand{\cI}{{\cal I}}
\newcommand{\cL}{{\cal L}}
\newcommand{\cM}{{\cal M}}
\newcommand{\cO}{{\cal O}}
\newcommand{\cP}{{\cal P}}
\newcommand{\fsu}{{\mathfrak s}{\mathfrak u}}
\newcommand{\cV}{{\cal V}}
\newcommand{\cW}{{\cal W}}
\newcommand{\Mbar}{\overline{\cM}}
\newcommand{\vac}{|0\rangle}
\newcommand{\lvac}{\langle 0|}

\DeclareMathOperator{\Aut}{Aut}
\DeclareMathOperator{\Ext}{Ext}
\DeclareMathOperator{\End}{End}
\DeclareMathOperator{\ext}{ext}
\DeclareMathOperator{\diag}{diag}
\DeclareMathOperator{\Span}{span}
\DeclareMathOperator{\id}{id}
\DeclareMathOperator{\val}{val}
\DeclareMathOperator{\sign}{sign}

\DeclareMathOperator{\res}{res}
\DeclareMathOperator{\tr}{tr}

\newtheorem{theorem}{Theorem}[section]
\newtheorem{Theorem}{Theorem}
\newtheorem{proposition}[theorem]{Proposition}
\newtheorem{lemma}[theorem]{Lemma}

\theoremstyle{remark}
\newtheorem{remark}{Remark}[section]

\theoremstyle{definition}

\begin{document}

\title{A proof of a Conjecture of
Mari\~{n}o-Vafa on Hodge Integrals}
\author{Chiu-Chu Melissa Liu}
\address{Department of Mathematics,
Harvard University, Cambridge, MA 02138, USA}
\email{ccliu@math.harvard.edu}
\author{Kefeng Liu}
\address{Center of Mathematical Sciences, Zhejiang University, Hangzhou, Zhejiang 310027, China; Department of Mathematics,
University of California at Los Angeles, Los Angeles, CA 90095-1555, USA}
\email{liu@cms.zju.edu.cn, liu@math.ucla.edu}
\author{Jian Zhou}
\address{Department of Mathematical Sciences\\
Tsinghua University\\ Beijing, 100084, China}
\email{jzhou@math.tsinghua.edu.cn}
\begin{abstract}
We prove a remarkable formula for Hodge integrals conjectured by
Mari\~{n}o and Vafa \cite{Mar-Vaf} based on large $N$ duality,
using functorial virtual localization on certain moduli spaces
of relative stable morphisms.
\end{abstract}

\maketitle

\section{Introduction}\label{introduction}

Let $\Mbar_{g,n}$ denote the Deligne-Mumford moduli stack of stable curves of genus
$g$ with $n$ marked points. Let $\pi:\Mbar_{g,n+1}\to \Mbar_{g,n}$
be the universal curve, and let $\omega_\pi$ be the relative dualizing sheaf.
The Hodge bundle
$$
\bE=\pi_*\omega_\pi
$$
is a rank $g$ vector bundle over $\Mbar_{g,n}$ whose
fiber over $[(C,x_1,\ldots,x_n)]\in\Mbar_{g,n}$ is
$H^0(C,\omega_C)$.
Let $s_i:\Mbar_{g,n}\to\Mbar_{g,n+1}$ denote the section of $\pi$ which
corresponds to the $i$-th marked point, and let
$$
\bL_i=s_i^*\omega_\pi
$$
be the line bundle over $\Mbar_{g,n}$ whose fiber over
$[(C,x_1,\ldots,x_n)]\in\Mbar_{g,n}$ is the cotangent line $T^*_{x_i} C$
at the $i$-th marked point $x_i$.
A Hodge integral is an integral of the form
$$\int_{\Mbar_{g, n}} \psi_1^{j_1}
\cdots \psi_n^{j_n}\lambda_1^{k_1} \cdots \lambda_g^{k_g}$$
where
$\psi_i=c_1(\bL_i)$ is the first Chern class of $\bL_i$, and
$\lambda_j=c_j(\bE)$ is the $j$-th Chern class of the Hodge bundle.

Hodge integrals arise in the calculations of Gromov-Witten invariants
by localization techniques \cite{Kon2, Gra-Pan}.
The explicit evaluation of Hodge integrals is a difficult problem.
The Hodge integrals involving only $\psi$ classes
can be computed recursively by Witten's conjecture \cite{Wit}
proved by Kontsevich \cite{Kon1}. Algorithms of computing Hodge
integrals are described in \cite{Fab}.

In \cite{Mar-Vaf}, M. Mari\~{n}o and C. Vafa obtained a closed
formula for a generating function of certain open Gromov-Witten
invariants, some of which has been reduced to Hodge integrals
by localization techniques which are not fully clarified
mathematically.
This leads to a conjectural formula of Hodge integrals by comparing
with the calculations in \cite{Kat-Liu}. To state this
formula, we introduce some notation.
Let
$$
\Lambda^{\vee}_g(u)=u^g-\lambda_1 u^{g-1} +\cdots+(-1)^g\lambda_g
$$
be the Chern polynomial of $\bE^\vee$, the dual of the
Hodge bundle. For a partition $\mu$ given by
$$
\mu_1\geq\mu_2\geq\cdots\geq\mu_{l(\mu)} >0,
$$
let $|\mu|=\sum_{i=1}^{l(\mu)}\mu_i$, and define
\begin{eqnarray*}
\cC_{g, \mu}(\tau)& = & - \frac{\sqrt{-1}^{|\mu|+l(\mu)}}{|\Aut(\mu)|}
 (\tau(\tau+1))^{l(\mu)-1}
\prod_{i=1}^{l(\mu)}\frac{ \prod_{a=1}^{\mu_i-1} (\mu_i \tau+a)}{(\mu_i-1)!} \\
&& \cdot \int_{\Mbar_{g, l(\mu)}}
\frac{\Lambda^{\vee}_g(1)\Lambda^{\vee}_g(-\tau-1)\Lambda_g^{\vee}(\tau)}
{\prod_{i=1}^{l(\mu)}(1- \mu_i \psi_i)}, \\
\cC_{\mu}(\lambda; \tau) & = & \sum_{g \geq 0} \lambda^{2g-2+l(\mu)}\cC_{g,
\mu}(\tau)
\end{eqnarray*}

Note that
$$
\int_{\Mbar_{0, l(\mu)}}
\frac{\Lambda^\vee_0(1)\Lambda^\vee_0(-\tau-1)\Lambda_0^\vee(\tau) }
{\prod_{i=1}^{l(\mu)} (1 -\mu_i \psi_i) }
=\int_{\Mbar_{0, l(\mu)}}
\frac{1}{\prod_{i=1}^{l(\mu)}(1 - \mu_i\psi_i)}
= |\mu|^{l(\mu)-3}
$$
for $l(\mu)\geq 3$, and we use this expression to extend the definition
to the case $l(\mu)<3$.

Introduce formal variables $p=(p_1,p_2,\ldots,p_n,\ldots)$, and define
$$
p_\mu=p_{\mu_1}\cdots p_{\mu_{l(\mu)} }
$$
for a partition $\mu$. Define
generating functions
\begin{eqnarray*}
\cC(\lambda; \tau; p) & = & \sum_{|\mu| \geq 1} \cC_{\mu}(\lambda;\tau)p_{\mu}, \\
\cC(\lambda; \tau; p)^{\bullet} & = & e^{\cC(\lambda; \tau; p)}.
\end{eqnarray*}

As pointed out in \cite{Mar-Vaf}, by comparing computations in
\cite{Mar-Vaf} with computations in \cite{Kat-Liu},
one obtains a conjectural formula for $\cC_{\mu}(\tau)$.
This formula is explicitly written down in \cite{Zho3}:
\begin{equation}\label{eqn:Mar-Vaf1}
\cC(\lambda; \tau; p)
=\sum_{n \geq 1} \frac{(-1)^{n-1}}{n}\sum_{\mu}
\left(\sum_{\cup_{i=1}^n \mu^i = \mu}
\prod_{i=1}^n \sum_{|\nu^i|=|\mu^i|} \frac{\chi_{\nu^i}(C(\mu^i))}{z_{\mu^i}}
e^{\sqrt{-1}(\tau+\frac{1}{2})\kappa_{\nu^i}\lambda/2} V_{\nu^i}(\lambda)
\right)p_\mu,
\end{equation}
\begin{equation}\label{eqn:Mar-Vaf2}
 \cC(\lambda;\tau; p)^{\bullet} = \sum_{|\mu| \geq 0}
\left(\sum_{|\nu|=|\mu|} \frac{\chi_{\nu}(C(\mu))}{z_{\mu}}
e^{\sqrt{-1}(\tau+\frac{1}{2})\kappa_{\nu}\lambda/2} V_{\nu}(\lambda)\right)
p_\mu,
\end{equation}
where
\begin{equation} \label{eqn:V}
\begin{split}
V_{\nu}(\lambda) = & \prod_{1 \leq a < b \leq l(\nu)}
\frac{\sin \left[(\nu_a - \nu_b + b - a)\lambda/2\right]}{\sin \left[(b-a)\lambda/2\right]} \\
& \cdot\frac{1}{\prod_{i=1}^{l(\nu)}\prod_{v=1}^{\nu_i} 2 \sin \left[(v-i+l(\nu))\lambda/2\right]}.
\end{split} \end{equation}
The right-hand side of (\ref{eqn:Mar-Vaf1}) is actually some truncated version of the
more general formula \cite[(5.6)]{Mar-Vaf} given by Mari\~{n}o and Vafa.

We now explain the notation on the right-hand sides of (\ref{eqn:Mar-Vaf1}) and
(\ref{eqn:Mar-Vaf2}). For a partition $\mu$,
$\chi_{\mu}$ denotes the character of the irreducible representation of $S_d$
indexed by $\mu$, where $d=|\mu| = \sum_{i=1}^{l(\mu)} \mu_i$.
The number $\kappa_{\mu}$ is defined by
$$\kappa_{\mu} = |\mu| + \sum_i (\mu_i^2 - 2i\mu_i).$$
For each positive integer $i$,
$$m_i(\mu) = |\{j: \mu_j = i\}|.$$
Denote by $C(\nu)$ the conjugacy class of $S_d$ corresponding to the partition $\nu$,
and by $\chi_{\mu}(C(\nu))$ the value of the character $\chi_{\mu}$ on the conjugacy class $C(\nu)$.
Finally,
$$z_{\mu} =  \prod_j m_j(\mu)!j^{m_j(\mu)}. $$

In this paper, we will call (\ref{eqn:Mar-Vaf1}) the Mari\~{n}o-Vafa formula.
The third author proved in \cite{Zho1} some special cases of the
Mari\~{n}o-Vafa formula and found some
applications \cite{Zho4, Zho5}.

We now describe our approach to the
Mari\~{n}o-Vafa formula (\ref{eqn:Mar-Vaf1}).
Denote the right-hand sides of (\ref{eqn:Mar-Vaf1})
and (\ref{eqn:Mar-Vaf2}) by $R(\lambda; \tau;p)$
and $R(\lambda;\tau;p)^{\bullet}$ respectively.
In \cite{Zho3}, the third author proved that the following two equivalent 
cut-and-join equations similar to the one satisfied by Hurwitz numbers
(see \cite{Gou-Jac-Vai}, \cite{Li-Zha-Zhe}, \cite[Section 15.2]{Ion-Par}):

\begin{Theorem}  \label{RCutJoin}
\begin{eqnarray}
&& \frac{\partial R}{\partial \tau}
= \frac{\sqrt{-1}\lambda}{2} \sum_{i, j\geq 1} \left(ijp_{i+j}\frac{\partial^2R}{\partial p_i\partial p_j}
+ ijp_{i+j}\frac{\partial R}{\partial p_i}\frac{\partial R}{\partial p_j}
+ (i+j)p_ip_j\frac{\partial R}{\partial p_{i+j}}\right), \label{eqn:CutJoin1} \\
&& \frac{\partial R^{\bullet}}{\partial \tau}
= \frac{\sqrt{-1}\lambda}{2} \sum_{i, j\geq 1} \left(ijp_{i+j}\frac{\partial^2R^{\bullet}}{\partial p_i\partial p_j}
+ (i+j)p_ip_j\frac{\partial R^{\bullet}}{\partial p_{i+j}}\right). \label{eqn:CutJoin2}
\end{eqnarray}
\end{Theorem}

Here is a crucial observation: One can rewrite (\ref{eqn:CutJoin2}) as a
sequence systems of ordinary equations, one for each positive integer $d$,
hence if $\cC(\lambda;\tau;p)^{\bullet}$ satisfies (\ref{eqn:CutJoin2}),
then it is determined by the initial value $\cC(\lambda; 0; p)^{\bullet}$.
To prove (\ref{eqn:Mar-Vaf1}) or (\ref{eqn:Mar-Vaf2}),
it suffices to prove the following two statements:
\begin{itemize}
\item[(a)] Equation (\ref{eqn:CutJoin1}) is satisfied by $\cC(\lambda; \tau; p)$.
\item[(b)] $\cC(\lambda; 0; p) = R(\lambda; 0; p)$.
\end{itemize}
Or equivalently,
\begin{itemize}
\item[(a)'] Equation (\ref{eqn:CutJoin2}) is satisfied by $\cC(\lambda; \tau; p)^\bullet$.
\item[(b)'] $\cC(\lambda; 0; p)^\bullet = R(\lambda; 0; p)^\bullet$.
\end{itemize}

The generating function $\cC(\lambda;0;p)$ of Hodge integrals has a closed
formula \cite[Theorem 2]{Fab-Pan1}.
This closed formula is shown to be equal to $R(\lambda;0;p)$ in \cite{Zho3}.
Therefore, the Mari\~{n}o-Vafa formula (\ref{eqn:Mar-Vaf1})
follows from the following theorem.

\begin{Theorem}\label{CCutJoin}
\begin{eqnarray}
&& \frac{\partial \cC}{\partial \tau}
= \frac{\sqrt{-1}\lambda}{2} \sum_{i, j\geq 1} \left(ijp_{i+j}
\frac{\partial^2\cC}{\partial p_i\partial p_j}
+ ijp_{i+j}\frac{\partial \cC}{\partial p_i}\frac{\partial\cC}{\partial p_j}
+ (i+j)p_ip_j\frac{\partial \cC}{\partial p_{i+j}}\right). \label{eqn:CutJoin}
\end{eqnarray}
\end{Theorem}

Our earlier paper \cite{LLZ} contains an essentially complete proof
of the Mari\~{n}o-Vafa formula based on the approach described above.
The purpose of this paper is to supply various computational details,
and to present some related results.
See \cite{Oko-Pan} for another approach to the Mari\~{n}o-Vafa formula.

The rest of this paper is organized as follows.
In Section \ref{initial}, we give a proof of the initial condition (b).
In Section \ref{combinatorics}, we give a proof of Theorem~\ref{RCutJoin}.
The materials in these two sections are already contained in \cite{Zho3} and some
 in \cite{LLZ}, and included here
for the convenience of the readers.
In Section \ref{moduli}, we recall the moduli spaces of relative stable
morphisms and obstruction bundles which will be used in the proof of Theorem
\ref{CCutJoin}.
In Section \ref{loci}, we use the graph notation to describe the torus fixed
points in the moduli spaces introduced in Section \ref{moduli}.
In Section \ref{geometry}, we prove Theorem \ref{CCutJoin} by functorial virtual
localization. Details  of virtual localization are
given in Appendix \ref{localization}.
In Section \ref{hurwitz}, we reproduce the cut-and-join equation of Hurwitz numbers
(proved in \cite{Gou-Jac-Vai}, \cite{Li-Zha-Zhe},
\cite[Section 15.2]{Ion-Par}) by virtual functorial localization.

\section{Initial Condition}\label{initial}

We recall some notations for partitions.
For a partition $\nu$, let $\nu'$ denote its transpose. Define
\begin{eqnarray} \label{eqn:neta}
n(\eta) = \sum_i (i-1)\eta_i
= \sum_i \begin{pmatrix} \eta'_i \\ 2\end{pmatrix}.
\end{eqnarray}
 The hook length of $\nu$ at the square $x$ located at the $i$-th row and $j$-th column
is defined to be:
$$h(x) = \nu_i + \nu_j'-i-j+1.$$
Then one has the following  two identities \cite[pp. 10 -11]{Mac}:
\begin{eqnarray}
&& \prod_{x \in \nu} (1 - t^{h(x)})
= \frac{\prod_{i=1}^{l(\nu)} \prod_{j=1}^{\nu_i - i + l(\nu)}( 1 - t^j)}
{\prod_{i < j} ( 1- t^{\nu_i - \nu_j - i +j})}, \label{eqn:hook1} \\
&& \sum_{x \in \nu} h(x) = n(\nu) + n(\nu') + |\nu|. \label{eqn:hook2}
\end{eqnarray}

We first obtain a simple expression for $V_{\nu}(\lambda)$.
\begin{proposition}\label{V}
$$
 V_{\nu}(\lambda) = \frac{1}{2^l\prod_{x \in \nu} \sin [h(x)\lambda/2]}.
$$
\end{proposition}

\begin{proof}
We rewrite the right-hand side of (\ref{eqn:hook1}) as follows.
\begin{eqnarray*}
&& RHS \\
& = &
\frac{\prod_{i=1}^{l(\nu)}\prod_{j=1}^{- i + l(\nu)}( 1 - t^j)}
    {\prod_{i < j} ( 1- t^{\nu_i - \nu_j - i + j})}
\cdot \prod_{i=1}^{l(\nu)}\prod_{j=1-i+l(\nu)}^{\nu_i - i + l(\nu)}( 1 - t^j) \\
& = &
\frac{\prod_{i<j} ( 1 - t^{j-i})}{\prod_{i < j} ( 1- t^{\nu_i - \nu_j - i +j})}
\cdot \prod_{i=1}^{l(\nu)}\prod_{j=1}^{\nu_i}( 1 - t^{j-i+l(\nu)}) \\
& = & t^{\left(\sum_{i<j} (j-i) - \sum_{i<j}(\nu_i-\nu_j-i+j) + \sum_{i=1}^{l(\nu)} \sum_{j=1}^{\nu_i} (j-i+l(\nu))\right)/2} \\
&& \cdot
\frac{\prod_{i<j} ( t^{-(j-i)/2} - t^{(j-i)/2})}{\prod_{i < j} ( t^{-(\nu_i-\nu_j-i+j)/2}- t^{(\nu_i - \nu_j - i +j)/2})}
\prod_{i=1}^{l(\nu)}\prod_{j=1}^{\nu_i}( t^{-(j-i+l(\nu))/2} - t^{(j-i+l(\nu))/2}).
\end{eqnarray*}
Now
\begin{eqnarray*}
&& \sum_{i<j} (j-i) - \sum_{i<j}(\nu_i-\nu_j-i+j) + \sum_{i=1}^{l(\nu)} \sum_{j=1}^{\nu_i} (j-i+l(\nu)) \\
& = & - \sum_{i <j}\nu_i + \sum_{i< j} \nu_j + \sum_{i=1}^{l(\nu)} \sum_{j=1}^{\nu_i} j
- \sum_{i=1}^{l(\nu)} \sum_{j=1}^{\nu_i} i + \sum_{i=1}^{l(\nu)} \sum_{j=1}^{\nu_i} l(\nu) \\
& = & - \sum_{i=1}^{l(\nu)} (l(\nu) - i) \nu_i + \sum_{j=1}^{l(\nu)} (j-1)\nu_j
+ \sum_{i=1}^{l(\nu)} \frac{\nu_i(\nu_i+1)}{2} -\sum_{i=1}^{l(\nu)} i\nu_i + |\nu|l(\nu)\\
& = & - |\nu|l(\nu) + \sum_{i=1}^{l(\nu)} i\nu_i + \sum_{j=1}^{l(\nu)} (j-1)\nu_j
+ \sum_{i=1}^{l(\nu)} \frac{\nu_i(\nu_i-1)}{2} + |\nu| - \sum_{i=1}^{l(\nu)} i\nu_i + |\nu|l(\nu) \\
& = &  \sum_{j=1}^{l(\nu)} (j-1)\nu_j + \sum_{i=1}^{l(\nu)} \frac{\nu_i(\nu_i-1)}{2} + |\nu| \\
& = & n(\nu) + n(\nu') + |\nu| \\
& = & \sum_{x \in \nu} h(x).
\end{eqnarray*}
Comparing with the left-hand side,
one then gets:
\begin{eqnarray*}
&& \prod_{x \in \nu} (t^{-h(x)/2} - t^{h(x)/2}) \\
& = & \frac{\prod_{i<j} ( t^{-(j-i)/2} - t^{(j-i)/2})}{\prod_{i < j} ( t^{-(\nu_i-\nu_j-i+j)/2}- t^{(\nu_i - \nu_j - i +j)/2})}
\prod_{i=1}^{l(\nu)}\prod_{j=1}^{\nu_i}( t^{-(j-i+l(\nu))/2} - t^{(j-i+l(\nu))/2}).
\end{eqnarray*}
The proof is completed by taking $t = e^{-\sqrt{-1}\lambda}$.
\end{proof}

\smallskip

We next compute the initial values $\cC(\lambda;0;p)$.
\begin{proposition}\label{sin}
$$
 \cC(\lambda; 0; p) = - \sum_{d > 0}  \frac{\sqrt{-1}^{d+1}p_d}{2d\sin (d\lambda/2)}.
$$
\end{proposition}

\begin{proof}
When $l(\mu) > 1$,
we clearly have
$$\cC_{\mu}(\lambda;0) = 0.$$
When $\mu = (d)$ we have
\begin{eqnarray*}
\cC_{(d)} (\lambda; 0)
& = & - \sum_{g \geq 0}
\lambda^{2g-1} \sqrt{-1}^{d+1} \frac{\prod_{a=1}^{d-1} (d\cdot 0 + a)}{(d-1)!}
\int_{\Mbar_{g, 1}} \frac{\Lambda_g^{\vee}(1)\Lambda_g^{\vee}(0)\Lambda_g^{\vee}(-1)}
{1- d\psi_1} \\
& = & - \frac{\sqrt{-1}^{d+1}}{d}
\left((d\lambda)^{-1}+ \sum_{g \geq 1}  (d\lambda)^{2g-1}
\int_{\Mbar_{g, 1}} \lambda_g \psi_1^{2g-2}\right) \\
& = &  - \frac{\sqrt{-1}^{d+1}}{d^2\lambda} \cdot \frac{d\lambda/2}{\sin (d\lambda/2)} \\
& = & - \frac{\sqrt{-1}^{d+1}}{2d\sin (d\lambda/2)}.
\end{eqnarray*}
In the second equality we have used the Mumford's relations \cite[5.4]{Mum}:
$$\Lambda_g^\vee(1)\Lambda_g^{\vee}(-1)=(-1)^g.$$
In the third equality we have used \cite[Theorem 2]{Fab-Pan1}.
This completes the proof.
\end{proof}

\medskip

\begin{proposition} \label{log}
We have the following identity:
$$
 \log\left(  \sum_{n \geq 0} \sum_{|\rho|=n}\sum_{|\eta|=|\rho|}
\frac{e^{\frac{1}{4}\kappa_{\rho}\sqrt{-1}\lambda}}{\prod_{e \in \rho} 2\sin (h(e)\lambda/2)}
\frac{\chi_{\rho}(\eta)}{z_{\eta}} p_{\eta}\right)
= - \sum_{d \geq 1} \frac{\sqrt{-1}^{d+1}p_d}{2d\sin(d\lambda/2)}. \label{eqn:Key}
$$
\end{proposition}

\smallskip

To prove Proposition~\ref{log}, we need the following two lemmata.
\begin{lemma}
Introducing formal variables $x_1, \dots, x_n, \dots$ such that
$$p_i(x_1, \dots, x_n, \dots) = x_1^i + \cdots + x_n^i + \cdots.$$
Then we have
\begin{eqnarray} \label{eqn:pg}
&& \sum_{n \geq 0} t^n\sum_{|\rho|=n}\sum_{|\eta|=|\rho|}
\frac{q^{n(\rho)}}{\prod_{e \in \rho} (1 - q^{h(e)})}
\frac{\chi_{\rho}(\eta)}{z_{\eta}} p_{\eta}
= \frac{1}{\prod_{i, j} (1 -tx_iq^{j-1})}.
\end{eqnarray}
\end{lemma}

\begin{proof}
Recall the following facts about Schur polynomials:
\begin{eqnarray}
&& s_{\rho}(x) = \sum_{|\eta|=|\rho|}
\frac{\chi_{\rho}(\eta)}{z_{\eta}} p_{\eta}(x), \label{eqn:s} \\
&& s_{\rho}(1, q, q^2, \dots)
= \frac{q^{n(\rho)}}{\prod_{e \in \rho} (1 - q^{h(e)})}, \\
&& \sum_{n \geq 0} t^n \sum_{|\rho|=n} s_{\rho}(x)s_{\rho}(y)
= \frac{1}{\prod_{i, j} (1 -tx_iy_j)}.
\end{eqnarray}
Combining the last two identities, one gets:
\begin{eqnarray*}
&& \sum_{n \geq 0} t^n \sum_{|\rho|=n}
\frac{q^{n(\rho)}}{\prod_{e \in \rho} (1 - q^{h(e)})}s_{\rho}(x)
= \frac{1}{\prod_{i,j} (1 - tx_iq^{j-1})}.
\end{eqnarray*}
The proof is completed by (\ref{eqn:s}).
\end{proof}

\begin{lemma}
For any partition $\rho$ we have
\begin{eqnarray}
&& \frac{1}{2} \sum_{e \in \rho} h(e) - n(\rho)
=  \frac{1}{4}\kappa_{\rho} + \frac{1}{2}|\rho|.
\end{eqnarray}
\end{lemma}

\begin{proof}
\begin{eqnarray*}
&& \frac{1}{2} \sum_{e \in \rho} h(e) - n(\rho) = \frac{1}{2} (n(\rho') - n(\rho) + |\rho|)  \\
& = & \frac{1}{2} (\sum_i \begin{pmatrix} \rho_i \\ 2 \end{pmatrix}
- \sum_i (i-1)\rho_i + |\rho|) \\
& = & \frac{1}{4} (\sum_i \rho_i(\rho_i - 1) - 2 \sum_i i\rho_i + 4|\rho|) \\
& = & \frac{1}{4}\kappa_{\rho} + \frac{1}{2}|\rho|.
\end{eqnarray*}
\end{proof}

\begin{proof} (of Proposition~\ref{log})\\
Let $q = e^{-\sqrt{-1}\lambda}$, and $t = \sqrt{-1}q^{1/2}$.
Then
\begin{eqnarray*}
&& \sum_{n\geq 0} t^n\sum_{|\rho|=n}\sum_{|\eta|=|\rho|}
\frac{q^{n(\rho)}}{\prod_{e \in \rho} (1 - q^{h(e)})}
\frac{\chi_{\rho}(\eta)}{z_{\eta}} p_{\eta} \\
& = & \sum_{n \geq 0} \sqrt{-1}^nq^{n/2} \sum_{|\rho|=n}\sum_{|\eta|=|\rho|}
\frac{q^{n(\rho)-\frac{1}{2}\sum_{e\in \rho} h(e)}}{\prod_{e \in \rho} (q^{-h(e)/2} - q^{h(e)/2})}
\frac{\chi_{\rho}(\eta)}{z_{\eta}} p_{\eta} \\
& = & \sum_{n \geq 0} \sqrt{-1}^nq^{n/2} \sum_{|\rho|=n}\sum_{|\eta|=|\rho|}
\frac{q^{-\frac{1}{4}\kappa_{\rho} - \frac{1}{2}n}}{\prod_{e \in \rho} (q^{-h(e)/2} - q^{h(e)/2})}
\frac{\chi_{\rho}(\eta)}{z_{\eta}} p_{\eta} \\
& = & \sum_{n \geq 0} \sum_{|\rho|=n}\sum_{|\eta|=|\rho|}
\frac{e^{\frac{1}{4}\kappa_{\rho}\sqrt{-1}\lambda}}{\prod_{e \in \rho} 2\sin (h(e)\lambda/2)}
\frac{\chi_{\rho}(\eta)}{z_{\eta}} p_{\eta}.
\end{eqnarray*}
Hence by (\ref{eqn:pg}),
\begin{eqnarray*}
&& \log\left(  \sum_{n \geq 0} \sum_{|\rho|=n}\sum_{|\eta|=|\rho|}
\frac{e^{\frac{1}{4}\kappa_{\rho}\sqrt{-1}\lambda}}{\prod_{e \in \rho} 2\sin (h(e)\lambda/2)}
\frac{\chi_{\rho}(\eta)}{z_{\eta}} p_{\eta}\right) \\
& = & \log\left(  \sum_{n\geq 0} t^n\sum_{|\rho|=n}\sum_{|\eta|=|\rho|}
\frac{q^{n(\rho)}}{\prod_{e \in \rho} (1 - q^{h(e)})}
\frac{\chi_{\rho}(\eta)}{z_{\eta}} p_{\eta}\right) \\
& = & \log  \frac{1}{\prod_{i, j} (1 -tx_iq^{j-1})}
= \sum_{i,j \geq 1} \sum_{d \geq 1} \frac{1}{d} t^d q^{d(j-1)} x_i^d \\
& = & \sum_{j \geq 1} \sum_{d \geq 1} \frac{1}{d} t^d q^{d(j-1)} p_d
= \sum_{d \geq 1} \frac{p_d}{d} \frac{t^d}{1 - q^d} \\
& = & - \sum_{d \geq 0} \frac{\sqrt{-1}^{d+1}p_d}{2d\sin(d\lambda/2)}.
\end{eqnarray*}
\end{proof}

By Proposition \ref{V} and Proposition \ref{log}, we have
\begin{eqnarray*}
R(\lambda;0;p) & = & \log \left(  \sum_{n \geq 0} \sum_{|\rho|=n}\sum_{|\eta|=|\rho|}
\frac{\chi_{\rho}(\eta)}{z_{\eta}} e^{\frac{1}{4}\kappa_{\rho}\sqrt{-1}\lambda}
V_{\rho}(\lambda) p_{\eta}\right)\\
&=& \log\left(  \sum_{n \geq 0} \sum_{|\rho|=n}\sum_{|\eta|=|\rho|}
\frac{e^{\frac{1}{4}\kappa_{\rho}\sqrt{-1}\lambda}}{\prod_{e \in \rho} 2\sin (h(e)\lambda/2)}
\frac{\chi_{\rho}(\eta)}{z_{\eta}} p_{\eta}\right) \\
&=& - \sum_{d \geq 1} \frac{\sqrt{-1}^{d+1}p_d}{2d\sin(d\lambda/2)}.
\end{eqnarray*}
By Proposition \ref{sin}, we have
$$
\cC(\lambda;0;p)=- \sum_{d \geq 1} \frac{\sqrt{-1}^{d+1}p_d}{2d\sin(d\lambda/2)}.
$$
So the initial condition (b) holds.

\section{Proof of Theorem \ref{RCutJoin} } \label{combinatorics}

Let $\mu, \eta$ be two partitions, both represented by Young diagrams.
We write $\eta \in J_{i,j}(\mu)$ and $\mu \in C_{i,j}(\eta)$
if $\eta$ is obtained from $\mu$ by removing a row of length $i$ and a row
of length $j$,
then adding a row of length $i+j$. It is easy to see that
\begin{eqnarray} \label{eqn:CJm}
&& \frac{m_{i+j}(\eta)}{\prod_k m_k(\eta)!} = \begin{cases}
\frac{m_i(\mu)m_j(\mu)}{\prod_k m_k(\mu)!}, & i \neq j, \\
\frac{m_i(\mu)(m_i(\mu)-1)}{\prod_k m_k(\mu)!}, & i =j.
\end{cases}
\end{eqnarray}

Recall
$$c_{\mu} = \sum_{g \in C_{\mu}} g$$
lies in the center of the group algebra $\bC S_d$, hence it acts
as a scalar $f_{\nu}(\mu)$ on any irreducible representation $R_{\nu}$.
In other words,
let $\rho: S_d \to \End R_{\nu}$ be the representation indexed by $\nu$,
then
$$\sum_{g \in C(\mu)} \rho_{\nu}(g) = f_{\nu}(\mu)\id.$$
We need the following interpretation of $\kappa_{\nu}$ in terms
of character.

\begin{lemma}
We have
$$\kappa_{\nu} = 2 f_{\nu}(C(2)),$$
where we use $C(2)$ to denote the class of transpositions.
\end{lemma}

\begin{proof}
By \cite[p. 118, Example 7]{Mac},
\begin{eqnarray*}
f_{\nu}(C(2)) & = & |C(2)| \frac{\chi_{\nu}(C(2))}{\dim R_{\nu}} = n(\nu') - n(\nu)\\
& = & \sum_{i=1}^{l(\nu)} \begin{pmatrix} \nu_i \\ 2\end{pmatrix}
- \sum_{i=1}^{l(\nu)} (i-1)\nu_i \\
& = & \frac{1}{2} \sum_{i=1}^{l(\nu)} (\nu_i^2 - 2i \nu_i + \nu_i) = \frac{1}{2} \kappa_{\nu}.
\end{eqnarray*}
In the above we have used (\ref{eqn:neta}).
\end{proof}

We need the following result:

\begin{lemma}
Suppose $h \in S_d$ has cycle type $\mu$.
The product $c_{(2)} \cdot h$ is a sum of elements of $S_d$ whose type is either a cut or a join of $\mu$.
More precisely,
there are $ijm_i(\mu)m_j(\mu)$ (when $i < j$) or
$i^2 m_i(\mu)(m_i(\mu)-1)/2$ (when $i=j$) elements obtained from $h$
by joining an $i$-cycle in $h$ to a $j$-cycle in $h$,
and there are $(i+j)m_{i+j}(\mu)$ (when $i < j$) or $im_{2i}(\mu)$ (when $i=j$) elements obtained from $h$ by
cutting an $(i+j)$-cycle into an $i$-cycle and a $j$-cycle.
\end{lemma}

\begin{proof}
Denote by $[s_1, \dots, s_k]$ a $k$-cycle.
Then
$$[s, t]\cdot [s, s_2, \dots, s_i, t, t_2, \dots, t_j] = [s, s_2, \dots, s_i][t, t_2, \dots, t_j],$$
i.e., an $i+j$-cycle is cut into an $i$-cycle and a $j$-cycle.
Conversely,
$$[s, t]\cdot [s, s_2, \dots, s_i][t, t_2, \dots, t_j] = [s, s_2, \dots, s_i, t, t_2, \dots, t_j],$$
i.e., an $i$-cycle and a $j$-cycle is joined to an $i+j$-cycle.
Hence, for a permutation $h$ of type $\mu$,
$c_{(2)} \cdot h$ is a sum of all elements obtained from $h$ by either a cut or a join.
Fix a pair of $i$-cycle and $j$-cycle of $h$,
there are $i \cdot j$ different ways to join them to an $(i+j)$-cycle.
Taking into the account of $m_i(\mu)$ choices of $i$-cycles, and $m_j(\mu)$ choices of $j$-cycles,
we get the number of different ways to obtain an element from $h$ by joining an $i$-cycle in $h$ to a $j$-cycle in $h$
is
$$\begin{cases}
ijm_i(\mu)m_j(\mu), & i < j \\
i^2 m_i(\mu)(m_i(\mu)-1)/2, & i = j. \end{cases}$$
Similarly,
fix an $(i+j)$-cycle of $h$,
there are $i+j$ different ways to cut it into an $i$-cycle and a disjoint $j$-cycle when $i < j$.
When $i=j$,
there are only $i$ different ways to cut it into two $i$-cycles.
And taking into account the number of $(i+j)$-cycles in $h$,
we get the number of different ways to obtain an element from $h$ by cutting an $(i+j)$-cycle
into an $i$-cycle and a $j$-cycle is
$$\begin{cases}
(i+j)m_{i+j}(\mu), & i < j, \\
im_{2i}(\mu), & i = j.\end{cases}$$
\end{proof}

For any $h \in S_d$ of cycle type $\mu$ we have
\begin{eqnarray*}
&& \sum_{\mu}  f_{\nu}(2) \frac{\chi_{\nu}(\mu)}{z_{\mu}}p_{\mu} \\
& = & \sum_{\mu} \tr [f_{\nu}(2) \id \cdot \rho_{\nu}(h)]
\cdot  \prod_i \frac{p_i^{m_i(\mu)}}{i^{m_i(\mu)} m_i(\mu)!} \\
& = & \sum_{\mu} \tr [\sum_{g \in C(2)} \rho_{\nu}(g) \cdot \rho_{\nu}(h)]
\cdot  \prod_i \frac{p_i^{m_i(\mu)}}{i^{m_i(\mu)} m_i(\mu)!} \\
& = &  \sum_{\mu}\tr \rho_{\nu}(\sum_{g \in C(2)} g\cdot h)
\cdot  \prod_i \frac{p_i^{m_i(\mu)}}{i^{m_i(\mu)} m_i(\mu)!} \\
& = & \sum_{\mu}\left( \sum_{i < j} \left(\sum_{\eta \in J_{i,j}(\mu)} ijm_i(\mu)m_j(\mu) \chi_{\nu}(\eta)
+ \sum_{\eta \in C_{i,j}(\mu)} (i+j)m_{i+j}(\mu) \chi_{\nu}(\eta)\right) \right. \\
&& \left. + \sum_i\left(\sum_{\eta \in J_{i,i}(\mu)} \frac{1}{2}i^2m_i(\mu)(m_i(\mu)-1) \chi_{\nu}(\eta)
+ \sum_{\eta \in C_{i,i}(\mu)} im_{2i}(\mu) \chi_{\nu}(\eta)\right)\right) \\
&& \cdot  \prod_i \frac{p_i^{m_i(\mu)}}{i^{m_i(\mu)} m_i(\mu)!} \\
& = & \frac{1}{2} \sum_{i,j}
\left((i+j)p_ip_j \frac{\partial}{\partial p_{i+j}}
+ ijp_{i+j}\frac{\partial}{\partial p_i}\frac{\partial}{\partial p_j}\right)
\sum_{\eta} \frac{\chi_{\nu}(\eta)}{z_{\eta}} p_{\eta}.
\end{eqnarray*}
In the last equality we have used (\ref{eqn:CJm}). It follows that
\begin{eqnarray*}
&& \frac{\partial R(\lambda;\tau;p)^{\bullet}}{\partial \tau} \\
& = & \frac{\sqrt{-1}\lambda}{2} \sum_{\mu, \nu}
\left( f_{\nu}(2) \frac{\chi_{\nu}(C(\mu))}{z_{\mu}} p_{\mu}\right)
e^{\sqrt{-1}(\tau+\frac{1}{2})\kappa_{\nu}\lambda/2}
V_{\nu}(\lambda)  \\
& = & \frac{\sqrt{-1}\lambda}{2}\left(ijp_{i+j}\frac{\partial}{\partial p_i}\frac{\partial}{\partial p_j}
+ (i+j)p_ip_j \frac{\partial}{\partial p_{i+j}}\right) \sum_{\eta} \frac{\chi_{\nu}(\eta)}{z_{\eta}} p_{\eta}
e^{\sqrt{-1}(\tau+\frac{1}{2})\kappa_{\nu}\lambda/2}
V_{\nu}(\lambda).
\end{eqnarray*}
This finishes the proof of Theorem 1.

\section{Moduli Spaces of Relative Stable Morphisms} \label{moduli}

In this section, we introduce the geometric objects
involved in the proof of Theorem \ref{CCutJoin}.

\subsection{Moduli space of relative morphisms}\label{sec:MP}
We first describe the moduli space of relative stable morphisms to
$\bP^1$ used in \cite{Li-Son}. The moduli spaces of algebraic relative stable
morphisms are constructed by J. Li \cite{Li1}.

We introduce some notations. For any nonnegative integer $m$,
let
$$
\bP^1[m]=\bP^1_{(0)}\cup\bP^1_{(1)}\cup\cdots\cup\bP^1_{(m)}
$$
be a chain of $m+1$ copies $\bP^1$, where $\bP^1_{(l)}$ is
glued to $\bP^1_{(l+1)}$ at $p_1^{(l)}$ for $0\leq l \leq m-1$.
The irreducible component $\bP^1_{(0)}$ will be referred to as the
root component, and the other irreducible components will be called the
bubble components. A point $p_1^{(m)}\neq p_1^{(m-1)}$ is fixed on
$\bP^1_{(m)}$. Denote by $\pi[m]: \bP^1[m] \to \bP^1$ the map which
is identity on the root component and contracts all the bubble components
to $p_1^{(0)}$. For $m>0$, let
$$
\bP^1(m)=\bP^1_{(1)}\cup\cdots\cup\bP^1_{(m)}
$$
denote the union of bubble components of $\bP^1[m]$.

Let $\mu$ be a partition of $d>0$. Let
$\Mbar_{g, 0}(\bP^1, \mu)$ be the moduli
space of morphisms
$$f: (C,x_1,\dots,x_{l(\mu)}) \to (\bP^1[m],p_1^{(m)}),$$
such that
\begin{enumerate}
\item $(C,x_1,\ldots,x_{l(\mu)})$ is a prestable curve
      of genus $g$ with $l(\mu)$ marked points.
      For convenience, we assume the marked points are \emph{unordered}.
\item $f^{-1}(p_1^{(m)})=\sum_{i=1}^{l(\mu)}\mu_i x_i$ as Cartier divisors,
      and $\deg(\pi[m]\circ f)=d$.
\item The preimage of each node in $\bP^1[m]$ consists of nodes of $C$.
      If $f(y)=p_1^{(l)}$ and $C_1$ and $C_2$ are two irreducible
      components of $C$ which intersect at $y$, then
      $f|_{C_1}$ and $f|_{C_2}$ have the same contact order to
      $p_1^{(l)}$ at $y$.
\item The automorphism group of $f$ is finite.
\end{enumerate}
Two such morphisms are isomorphic if they differ by an isomorphism
of the domain and an automorphism of the pointed curve
$(\bP^1(m),p_1^{(0)},p_1^{(m)})$.
In particular, this defines the automorphism group in the stability
condition (4) above.

In \cite{Li1, Li2}, J. Li showed that $\Mbar_{g,0}(\bP^1,\mu)$ is a
separated, proper Deligne-Mumford stack with
a perfect obstruction theory of virtual
dimension
$$r=2g-2 + d + l(\mu),$$
so it has a virtual fundamental class of degree $r$.

\subsection{Torus action}
Consider the $\bC^*$-action
$$t \cdot [z^0:z^1] = [tz^0: z^1]$$
on $\bP^1$. It has two fixed points
$p_0= [0:1]$ and $p_1=[1:0]$.
This induces an action on $\bP^1[m]$ by the action on the root component
induced by the isomorphism to $\bP^1$, and the trivial actions on the
bubble components. This in turn induces an action on
$\Mbar_{g, 0}(\bP^1, \mu)$.

\subsection{The branch morphism}\label{branch}
There is a branch morphism
$$
\mathrm{Br}: \Mbar_{g,0}(\bP^1,\mu) \to
\mathrm{Sym}^r\bP^1\cong \bP^r.
$$
Note that $\bP^r$ can be identified with $\bP(H^0(\bP^1,\cO(r))$, and
the isomorphism
$$\bP(H^0(\bP^1,\cO(r))\cong\mathrm{Sym}^r\bP^1$$
is given by $[s]\mapsto \mathrm{div}(s)$.
The $\bC^*$-action on $\bP^1$ induces a
$\bC^*$-action on $H^0(\bP^1,\cO(r))$ by
$$
t\cdot (z^0)^k (z^1)^{r-k} = t^{-k}(z^0)^k (z^1)^{r-k}.
$$
So $\bC^*$ acts on $\bP^r$ by
$$
t\cdot [a_0:a_1:\cdots: a_r]= {[a_0:t^{-1}a_1:\cdots: t^{-r} a_r]},
$$
where $[a_0: a_1 :\ldots : a_r]$ corresponds to
$\sum_{k=0}^r a_k (z^0)^k (z^1)^{r-k}\in H^0(\bP^1,\cO(r))$.
With this action, the branch morphism is $\bC^*$-equivariant.
See \cite{Fan-Pan, Gra-Vak1} for more details.

The $\bC^*$-action on $\bP^r$ has $r+1$ fixed points
$p_0,\ldots,p_r$, where $p_k\in \bP^r$ corresponds to the complex line
$\bC (z^0)^k(z^1)^{r-k}\subset H^0(\bP^1,\cO(r))$.

\subsection{The Obstruction Bundle}\label{obstruction}

In \cite{Li-Son}, J. Li and Y. Song constructed
an obstruction bundle over the stratum where the target
is $\bP^1[0]=\bP^1$, and proposed an extension over the
entire $\Mbar_{g,0}(\bP^1,\mu)$.
Here we use a different extension which is defined in
\cite[Section 3]{Bry-Pan}.

Let
$$
\pi:\mathcal{U}_{g,\mu}\to\Mbar_{g,0}(\bP^1,\mu)
$$
be the universal domain curve, and let
$$
P:\mathcal{T}_{g,\mu}\to\Mbar_{g,0}(\bP^1,\mu)
$$
be the universal target.
There is an evaluation map
$$
F:\mathcal{U}_{g,\mu}\to \mathcal{T}_{g,\mu}
$$
and a contraction map
$$
\tilde{\pi}:\mathcal{T}_{g,\mu}\to \bP^1.
$$
Let $\mathcal{D}_{g,\mu}\subset \mathcal{U}_{g,\mu}$ be the
divisor corresponding to the $l(\mu)$ marked points.
Define
\begin{eqnarray*}
V_D&=&R^1\pi_*(\cO_{\mathcal{U}_{g,\mu}}(-\mathcal{D}_{g,\mu}) )\\
V_{D_d}&=&R^1\pi_* \tilde{F}^*\cO_{\bP^1}(-1),
\end{eqnarray*}
where $\tilde{F}=\tilde\pi\circ F:\mathcal{U}_{g,\mu}\to \bP^1$.
The fibers of $V_D$ and $V_{D_d}$ at
$$
\left[ f:(C,x_1,\ldots,x_{l(\mu)})\to \bP^1[m]\ \right]\in \Mbar_{g,0}(\bP^1,\mu)
$$
are $H^1(C, \cO_C(-D))$ and $H^1(C, \tilde{f}^*\cO_{\bP^1}(-1))$,
respectively, where $D=x_1+\ldots+x_{l(\mu)}$, and $\tilde{f}=\pi[m]\circ f$.
Note that
$$H^0(C, \cO_C(-D))= H^0(C, \tilde{f}^*\cO_{\bP^1}(-1))=0,$$
so  $V_D$ and $V_{D_d}$ are vector bundles of ranks
$l(\mu)+g-1$ and $d+g-1$, respectively.
The obstruction bundle
$$
V=V_D\oplus V_{D_d}
$$
is a vector bundle of rank $r=2g-2+d+l(\mu)$.

We lift the $\bC^*$-action on $\Mbar_{g,0}(\bP^1,\mu)$ to
$V_D$ and $V_{D_d}$ as follows. The action on
$V_{D_d}$ comes from an action on $\cO_{\bP^1}(-1)\to\bP^1$
with weights $-\tau-1$ and $-\tau$ at the two fixed points $p_0$ and
$p_1$, respectively, where $\tau\in\bZ$. The fiber of
$V_D$ does not depend on the map $f$, so the fibers over
two points in the same orbit of the $\bC^*$-action can be
canonically identified. The action of $\lambda\in\bC^*$ on $V_D$
is multiplication by $\lambda^\tau$.

\section{Fixed Points of Torus Action}\label{loci}

\subsection{Graph notation} \label{graph}
Similar to the case of $\Mbar_{g, 0}(\bP^1, d)$,
the connected components of the $\bC^*$ fixed points set
$\Mbar_{g, 0}(\bP^1, \mu)^{\bC^*}$ are parameterized by labeled graphs.

Given a morphism
$$ f:(C,x_1,\ldots,x_{l(\mu)})\to \bP^1[m] $$
which represents a fixed point of the $\bC^*$-action
on $\Mbar_{g,0}(\bP^1,\mu)$, let
$$ \tilde{f}=\pi[m]\circ f: C\to \bP^1. $$
The restriction of $\tilde{f}$ to an irreducible
component of $C$ is either a constant map to one of the $\bC^*$
fixed points $p_0, p_1$ or a cover of $\bP^1$ which is fully
ramified over $p_0$ and $p_1$. We associate a labeled graph $\Gamma$ to
the $\bC^*$ fixed point
$$\left[ f:(C,x_1,\ldots,x_{l(\mu)})\to \bP^1[m]\ \right]$$ as follows:
\begin{enumerate}
\item We assign a vertex $v$ to each connected
component $C_v$ of $\tilde{f}^{-1}(\{p_0,p_1\})$, a label
$i(v)=i$ if $\tilde{f}(C_v)=p_i$, where $i=0,1$, and a label $g(v)$
which is the arithmetic genus of $C_v$ (We define $g(v)=0$ if $C_v$ is a
point). Denote by $V(\Gamma)^{(i)}$ the set of vertices with $i(v)=i$,
where $i=0,1$. Then the set $V(\Gamma)$ of vertices of the graph $\Gamma$
is a disjoint union of $V(\Gamma)^{(0)}$ and $V(\Gamma)^{(1)}$.

\item We assign an edge $e$ to each rational irreducible component $C_e$ of $C$
such that $\tilde{f}|_{C_e}$ is not a constant map.
Let $d(e)$ be the degree of $\tilde{f}|_{C_e}$.
Then $\tilde{f}|_{C_e}$ is fully ramified over $p_0$ and $p_1$.
Let $E(\Gamma)$ denote the set of edges of $\Gamma$.

\item The set of flags of $\Gamma$ is given by
$$F(\Gamma)=\{(v,e):v\in V(\Gamma), e\in E(\Gamma),
C_v\cap C_e\neq \emptyset \}.$$

\item For each $v\in V(\Gamma)$, define
$$d(v)=\sum_{(v,e)\in F(\Gamma)}d(e),$$
and let $\nu(v)$ be the partition of $d(v)$
determined by $\{d(e): (v,e)\in F(\Gamma)\}$.
When the target is $\bP^1[m]$, where $m>0$,
we assign an additional label for each
$v\in V(\Gamma)^{(1)}$:
let $\mu(v)$ be the partition of
$d(v)$ determined by the ramification of
$f|_{C_v}:C_v\to\bP^1(m)$ over $p_1^{(m)}$.
\end{enumerate}
Note that for $v\in V(\Gamma)^{(1)}$,
$\nu(v)$ coincides with the partition of $d(v)$
determined by the ramification of $f|_{C_v}:C_v\to \bP^1(m)$
over $p_1^{(0)}$.

\subsection{Fixed points}
Let $G_{g,0}(\bP^1,\mu)$ be the set of
all the graphs associated to the $\bC^*$ fixed
points in $\Mbar_{g,0}(\bP^1,\mu)$, as described
in Section \ref{graph}.
In this section, we describe the set of fixed
points associated to a given graph $\Gamma\in G_{g,0}(\bP^1,\mu)$.

\subsubsection{The target is $\bP^1$}
Any $\bC^*$ fixed point in $\Mbar_{g,0}(\bP^1,\mu)$
which is represented by a morphism to $\bP^1$ is
associated to the graph $\Gamma^0$, where
$$
V(\Gamma^0)^{(0)}=\{v_0\},\ \ \
V(\Gamma^0)^{(1)}=\{v_1,\ldots,v_{l(\mu)}\},\ \ \
E(\Gamma^0)=\{ e_1,\ldots,e_{l(\mu)} \},
$$
and
$$
g(v)=g,\ \ \
g(v_i)=0,\ \ \
d(e_i)=\mu_i
$$
for $i=1,\ldots,l(\mu)$.
The two end points of the edge $e_i$ are $v_0$ and $v_i$.
Let $\Aut(\mu)$ denote the automorphism group of the partition
$\mu$ of $d$. Any morphism associated to the graph $\Gamma^0$ has automorphism
group $A_{\Gamma^0}$, where
$$
1\to \prod_{i=1}^{l(\mu)}\bZ_{\mu_i}
\to A_{\Gamma^0} \to \Aut(\mu) \to 1.
$$

Let
$$
\Mbar_{\Gamma^0}= \begin{cases}
\{\mathrm{point}\}, & (g,l(\mu))=(0,1),(0,2),\\
\Mbar_{g,l(\mu)}, & (g,l(\mu))\neq (0,1),(0,2).
\end{cases}
$$
There is a morphism
$$
i_{\Gamma^0}:\Mbar_{\Gamma^0}\to
\Mbar_{g,0}(\bP^1,\mu)
$$
whose image is the fixed locus
$F_{\Gamma^0}$ associated to
$\Gamma_0$. The morphism $i_{\Gamma_0}$
induces an isomorphism
$$
\Mbar_{\Gamma^0}/A_{\Gamma^0}\cong F_{\Gamma^0}.
$$
The dimension of $F_{\Gamma^0}$ is
$$
d_{\Gamma^0}=
\left\{\begin{array}{ll}
0,& (g,l(\mu))=(0,1),(0,2),\\
3g-3+l(\mu),& (g,l(\mu))\neq (0,1),(0,2).
\end{array}\right.
$$

\subsubsection{The target is $\bP^1[m]$, $m>0$}

Let $\Gamma\in G_{g,0}(\bP^1,\mu)$ be a graph associated
to a $\bC^*$ fixed point represented by
some morphism to $\bP^1[m]$, $m>0$.
Let $\Mbar_{\Gamma}^{(1)}$ denote the moduli space of morphisms
$$
\hat{f}:\hat{C}\to \bP^1(m)
$$
such that
\begin{itemize}
\item[(a)] $\hat{C}$ is the disjoint union of $\{C_v:v\in V(\Gamma)^{(1)}\}$.
\item[(b)] $(C_v, x_{v,1},\ldots,x_{v,\l(\mu(v))},
             y_{v,1},\ldots,y_{v,\l(\nu(v))})$ is a prestable
      curve of genus $g(v)$ with $l(\mu(v))+l(\nu(v))$  marked points.
      Here the marked points are \emph{ordered}.
\item[(c)] As Cartier divisors,
$$
  (\hat{f}|_{C_v})^{-1}(p_1^{(0)})=\sum_{i=1}^{l(\mu(v))}\mu(v)_i\, x_{v,i},\ \ \
  (\hat{f}|_{C_v})^{-1}(p_1^{(m)})=
  \sum_{j=1}^{l(\nu(v))}\nu(v)_i\, y_{v,j}.
$$
The morphism $(\hat{f}|_{C_v})^{-1}(B)\to B$ is of degree $d(v)$ for each
irreducible component $B$ of $\bP^1(m)$.
\item[(d)] The automorphism group of $\hat{f}$ is finite.
\end{itemize}
Two such morphisms are isomorphic if they differ by an isomorphism
of the domain and an automorphism of the pointed curve
$(\bP^1(m),p_1^{(0)},p_1^{(m)})$, which is an element of $(\bC^*)^m$.
In particular, this defines the automorphism group in the stability
condition (d) above.

The moduli space $\Mbar_{\Gamma}^{(1)}$ is a
variant of J. Li's moduli spaces of stable relative morphisms
\cite{Li1, Li2}.
It is a separated, proper Deligne-Mumford stack with a perfect
obstruction theory.

Given
$$
f:(C,x_1,\ldots,x_{l(\mu)})\to \bP^1[m]
$$
associated to the graph $\Gamma$, let $\hat{C}$
be the disjoint union of $\{C_v:v\in V(\Gamma)^{(1)}\}$.
Let $\hat{f}$ be the restriction of $f$ to $\hat{C}$.
Then
$$
\hat{f}:\hat{C}\to \bP^1(m)
$$
represents a point in $\Mbar_\Gamma^{(1)}$.

Define
\begin{eqnarray*}
r_0(v)=&2g(v)-2 + \val(v), & v\in V(\Gamma)^{(0)},\\
r_1(v)=&2g(v)-2 + l(\mu(v))+l(\nu(v)), & v\in V(\Gamma)^{(1)},
\end{eqnarray*}
\begin{eqnarray*}
V^I(\Gamma)^{(0)}&=&\{v\in V(\Gamma)^{(0)}:r_0(v)=-1\},\\
V^{II}(\Gamma)^{(0)}&=&\{v\in V(\Gamma)^{(0)}:r_0(v)=0\},\\
V^S(\Gamma)^{(0)}&=&\{v\in V(\Gamma)^{(0)}:r_0(v)> 0\},\\
V^{II}(\Gamma)^{(1)}&=&\{v\in V(\Gamma)^{(1)}:r_1(v)=0\},\\
V^S(\Gamma)^{(1)}&=&\{v\in V(\Gamma)^{(1)}:r_1(v)>0\}.
\end{eqnarray*}
Note that $V^S(\Gamma)^{(1)}\neq\emptyset$ by the stability
condition (d).

Let $\Aut(\Gamma)$ denote the automorphism of the labeled
graph $\Gamma$.
The automorphism group of any morphism associated to
the graph $\Gamma$ is $A_\Gamma$, where
\[
1\to \prod_{e\in E(\Gamma)}\bZ_{d(e)} \to A_\Gamma\to \Aut(\Gamma)\to 1.
\]

Let $\Mbar_\Gamma=\Mbar_\Gamma^{(0)}\times \Mbar_\Gamma^{(1)}$, where
$$
\Mbar_{\Gamma}^{(0)} = \prod_{v \in V^S(\Gamma)^{(0)}}
\Mbar_{g(v),\val(v)}.
$$
There is a morphism
$$
i_{\Gamma}: \Mbar_{\Gamma}\to \Mbar_{g, 0}(\bP^1, \mu)
$$
whose image is the fixed locus $F_\Gamma$ associated to the
graph $\Gamma$. The morphism $i_\Gamma$ induces an isomorphism
$\Mbar_{\Gamma}/A_{\Gamma}\cong F_\Gamma$.

The dimension of $\Mbar_\Gamma^{(0)}$ is given by
$$
d_\Gamma^{(0)}=\sum_{v \in V^S(\Gamma)^{(0)}} (3g(v) - 3 + \val(v)),
$$
and the virtual dimension of  $\Mbar_\Gamma^{(1)}$ is given by
$$
d_\Gamma^{(1)}=\left(\sum_{v\in V(\Gamma)^{(1)} }r_1(v)\right) -1.
$$
So the virtual dimension of $F_\Gamma$ is given by
\begin{eqnarray*}
d_\Gamma
&=& d_\Gamma^{(0)}+ d_\Gamma^{(1)}\\
&=& \sum_{v \in V^S(\Gamma)^{(0)}} (3g(v) - 3 + \val(v))
   +\sum_{v\in V(\Gamma)^{(1)} }r_1(v)\, -1 \\
& = & \sum_{v \in V(\Gamma)^{(0)}} (3g(v) - 3 + \val(v))
+ |V^{II}(\Gamma)^{(0)}|+ 2|V^I(\Gamma)^{(0)}| \\
& & +\sum_{v \in V(\Gamma)^{(1)}} (2g(v)-2+l(\mu(v)) + l(\nu(v))) -1\\
& = & 3\sum_{v \in V(\Gamma)^{(0)}}g(v) - 3|V(\Gamma)^{(0)}| + |E(\Gamma)|
      +|V^{II}(\Gamma)^{(0)}|+ 2|V^I(\Gamma)^{(0)}| \\
&& +2\sum_{v \in V(\Gamma)^{(1)}}g(v) -2|V(\Gamma)^{(1)}|+l(\mu)+
|E(\Gamma)|-1\\
& = & 2\left(\sum_{v\in V(\Gamma)} g(v) -|V(\Gamma)|+|E(\Gamma)| +1 \right) -3 +l(\mu)\\
&   & +\sum_{v\in V(\Gamma)^{(0)}}(g(v)-1)+ |V^{II}(\Gamma)^{0}|
      +2|V^I(\Gamma)^{0}|\\
&=& 2g-3+l(\mu)+\sum_{v\in V^S(\Gamma)^{(0)}}(g(v)-1)+|V^I(\Gamma)^{0}|
\end{eqnarray*}
The last equality comes from the following identity:
$$
g=\sum_{v\in V(\Gamma)} g(v) + b_1(\Gamma)
 =\sum_{v\in V(\Gamma)} g(v) - |V(\Gamma)| + |E(\Gamma)| + 1,
$$
where $b_1(\Gamma)$ is the first betti number of the graph $\Gamma$.

\section{Proof of Theorem \ref{CCutJoin}} \label{geometry}
\subsection{Functorial localization}
Let $T=\bC^*$. We have seen in Section \ref{branch} that
the branch morphism
$$\mathrm{Br}:\Mbar_{g,0}(\bP^1,\mu)\to\bP^r $$
is $T$-equivariant. We will compute
$$
\mathrm{Br}_* e_T(V)=\sum_{l=0}^r a_l(\tau) H^l u^{r-l}.
$$
by virtual functorial localization \cite{LLY},
where $H\in H^2(\bP^r;\bZ)$ is the hyperplane class, and $a_l(\tau)$ is
a polynomial in $\tau$. Recall that $\tau\in\bZ$ parametrizes
torus actions on the obstruction bundle, as described in Section \ref{obstruction}.

Let $p_0,\ldots, p_r\in \bP^r$ be the torus fixed points
defined as in Section \ref{branch},
and let $f_k:p_k\to\bP^r$ be the inclusion.
From the torus action on $\bP^r$ described in Section \ref{branch}, one gets
$$
\frac{f_k^*\mathrm{Br}_* e_T(V)}{e_T(T_{p_k}\bP^r)}=
\frac{F(\tau,k)}{(-1)^{r-k} k!(r-k)!},
$$
where
$$F(\tau,x)=\sum_{l=0}^r a_l(\tau)x^l.$$

By functorial localization, we have
$$
\int_{p_k}\frac{f_k^* \mathrm{Br}_* e_T(V)}{e_T(T_{p_k}\bP^r)}
=\sum_{F_\Gamma\subset\mathrm{Br}^{-1}(p_k)}
\frac{1}{|A_\Gamma|}\int_{[\Mbar_\Gamma]^{\mathrm{vir} }}
\frac{i_\Gamma^* e_T(V)}{e_T(N_\Gamma^{\mathrm{vir}}) }
$$
for $k=0,\ldots,r$, where $N_\Gamma^{\mathrm{vir}}\to \Mbar_\Gamma$ is the
pull-back of the virtual normal bundle of $F_\Gamma$ in $\Mbar_{g,0}(\bP^1,\mu)$.
Note that $\mathrm{Br}(F_{\Gamma^0})=p_r$, and
$$
\mathrm{Br}(F_\Gamma)=p_{r-d_\Gamma^{(1)}-1}
$$
for $\Gamma\in G_{g,0}(\bP^1,\mu)$, $\Gamma\neq\Gamma^0$.
Recall that  $d_\Gamma^{(1)}$ is the virtual dimension of
$\Mbar_\Gamma^{(1)}$, and $0\leq d_\Gamma^{(1)}\leq r-1$.
So
\begin{eqnarray*}
F(\tau,x)&=&\sum_{k=0}^r\frac{F(\tau,k)}{(-1)^{r-k} k!(r-k)!}
x(x-1)\cdots(x-k+1)(x-k-1)\cdots (x-r)\\
&=&\sum_{k=0}^r I_{g,\mu}^{r-k}(\tau)
x(x-1)\cdots(x-k+1)(x-k-1)\cdots (x-r),
\end{eqnarray*}
where
$$
I_{g,\mu}^0(\tau)=\frac{1}{|A_{\Gamma^0}|}\int_{\Mbar_{\Gamma^0}}
\frac{i_{\Gamma^0}^* e_T(V)}{e_T(N_{\Gamma^0}^{\mathrm{vir}}) },
$$
and
$$
I_{g,\mu}^k(\tau)=\sum_{\Gamma\in G_{g,0}(\bP^1,\mu),
                     \Gamma\neq \Gamma^0, d_\Gamma +1=k}
\frac{1}{|A_\Gamma|}\int_{[\Mbar_{\Gamma}]^{\mathrm{vir} }}
\frac{i_\Gamma^* e_T(V)}{e_T(N_{\Gamma}^{\mathrm{vir}}) }
$$
for $k=1,\ldots,r$.

\subsection{Contribution from each graph}

\subsubsection{The target is $\bP^1$}
Consider the graph $\Gamma^0\in G_{g,0}(\bP^1,\mu)$.
We first consider the stable case, i.e., $(g,l(\mu))\neq (0,1),(0,2)$.
Let $d=|\mu|$ as before.
Using the Feynman rules derived in Appendix \ref{localization}, we obtain
\begin{eqnarray*}
I_{g,\mu}^0(\tau)&=&
\frac{1}{|A_{\Gamma^0}|}\int_{\Mbar_{\Gamma^0}}
\frac{i_{\Gamma^0}^* e_T(V)}{e_T(N_{\Gamma^0}^{\mathrm{vir}}) }\\
& = &\frac{(-1)^{d-1} }{|\mathrm{Aut}(\mu)|}
      (\tau(\tau+1))^{l(\mu)-1}
      \left(\prod_{i=1}^{l(\mu)}\frac{\prod_{a=1}^{\mu_i-1}(\mu_i\tau+a)}{(\mu_i-1)!}
     \right) \\
&   & \cdot \int_{\Mbar_{g,l(\mu)} }
    \frac{\Lambda_g^{\vee}(u) \Lambda_{g}^{\vee}(\tau u)
    \Lambda_{g}^{\vee}((-\tau-1)u) u^{2l(\mu)-3} }{\prod_{i=1}^{l(\mu)}
    (u-\mu_i\psi_i)}\\
& = & \frac{(-1)^{d-1} }{|\mathrm{Aut}(\mu)|} (\tau(\tau+1))^{l(\mu)-1}
 \left(\prod_{i=1}^{l(\mu)}
   \frac{\prod_{a=1}^{\mu_i-1}(\mu_i\tau+a)}{(\mu_i-1)!}\right) \\
&   & \cdot \int_{ \Mbar_{g,l(\mu)} }
    \frac{\Lambda_g^{\vee}(1) \Lambda_{g}^{\vee}(\tau)
    \Lambda_{g}^{\vee}(-\tau-1)}{\prod_{i=1}^{l(\mu)}
    (1-\mu_i\psi_i)}\\
& = & \sqrt{-1}^{d-l(\mu)}\cC_{g,\mu}(\tau)
\end{eqnarray*}

In particular,
\begin{eqnarray*}
I^0_{0,\mu}(\tau)
&=&\frac{(-1)^{d-1}}{|\mathrm{Aut}(\mu)|}(\tau(\tau+1))^{l(\mu)-1}
        \left(
     \prod_{i=1}^{l(\mu)}\frac{\prod_{a=1}^{\mu_i-1}(\mu_i\tau+a)}{(\mu_i-1)!}
     \right) \int_{ \Mbar_{0,l(\mu)} }
     \frac{1}{\prod_{i=1}^{l(\mu)}(1-\mu_i\psi_i)}\\
&=& \frac{ (-1)^{d-1} }{|\mathrm{Aut}(\mu)|}(\tau(\tau+1))^{l(\mu)-1}
  \left(\prod_{i=1}^{l(\mu)}\frac{\prod_{a=1}^{\mu_i-1}(\mu_i\tau+a)}{(\mu_i-1)!}
  \right) d^{l(\mu)-3}.
\end{eqnarray*}

The contribution from $\Gamma^0$ for the two unstable cases is also given by the
above formula:
\begin{eqnarray*}
      I^0_{0,(d)}(\tau)
& = & (-1)^{d-1}\left( \frac{\prod_{a=1}^{d-1}(d\tau+a)}{(d-1)!}
       \right)\frac{1}{d^2}.\\
       I^0_{0,(\mu_1,\mu_2)}(\tau)
& = & \frac{(-1)^{d-1}}{|\mathrm{Aut}((\mu_1,\mu_2))|}\tau(\tau+1)\left(
      \prod_{i=1}^{2}\frac{\prod_{a=1}^{\mu_i-1}(\mu_i\tau+a)}{(\mu_i-1)!}
      \right) \frac{1}{d}
\end{eqnarray*}

Note that $I^0_{g,\mu}(\tau)$ is a degree $r=2g-2+d+l(\mu)$ polynomial in
$\tau$ with rational coefficients, and
$$
I^0_{g,\mu}(-\tau-1)=(-1)^{d-l(\mu)}I^0_{g,\mu}(\tau).
$$

\subsubsection{The target is $\bP^1[m]$, $m>0$}

Consider $\Gamma\in G_{g,0}(\bP^1,\mu)$, $\Gamma\neq\Gamma_0$.
Using the Feynman rules derived in Appendix \ref{localization}, we obtain
$$
I^\Gamma = \frac{1}{|A_\Gamma|}\int_{[\Mbar_{\Gamma}]^{\mathrm{vir} }}
\frac{i_{\Gamma}^* e_T(V)}{e_T(N_{\Gamma}^{\mathrm{vir}}) }\\
=\frac{1}{|A_\Gamma|}\int_{[\Mbar_\Gamma]^{\mathrm{vir}} }
    \frac{\prod_{v\in V(\Gamma)} B_v }{-u-\psi^t}
$$
where
$$B_v=
\begin{cases}
A_v A_v^V\prod_{(v,e)\in F(\Gamma)}( A_e A_e^V),
& v\in V(\Gamma)^{(0)},\\
A_v A_v^V, &v\in V(\Gamma)^{(1)}.
\end{cases}
$$
More explicitly, in the notation of Appendix \ref{localization},
\begin{eqnarray*}
B_v
& = & \left(\prod_{(v,e)\in F(\Gamma)}d(e)\right)(-1)^{d(v)-1}
      (\tau(\tau+1))^{\val(v)-1}\left(
    \prod_{(v,e)\in F(\Gamma)}\frac{\prod_{a=1}^{d(e)-1}(d(e)\tau+a)}{(d(e)-1)!}
    \right)\\
&   & \cdot\frac{\Lambda^\vee_{g(v)}(u)\Lambda^\vee_{g(v)}(\tau u)
      \Lambda^\vee_{g(v)}(-(\tau +1)u)u^{2\val(v)-3} }{\prod_{(v,e)\in F(\Gamma)}
      \left(u-d(e)\psi_{(v,e)}\right)},
      \ \ \ v\in V^S(\Gamma)^{(0)},\\
&   & d(v)(-1)^{d(v)-1}\left(\frac{\prod_{a=1}^{d(v)-1}(d(v)\tau+a)}{(d(e)-1)!}\right)
       \frac{1}{d(v)^2},
      \ \ \ v\in V^I(\Gamma)^{(0)},\\
&   & d(e_1)d(e_2)(-1)^{d(v)-1}
       (\tau(\tau+1)) \left(
       \prod_{i=1}^2\frac{\prod_{a=1}^{d(e_i)-1}(d(e_i)\tau +a)}{(d(e_i)-1)!}
       \right)\frac{1}{d(v)},\\
&   &  \ \ \   v\in V^{II}(\Gamma)^{(0)}, (v,e_1),(v,e_2)\in F(\Gamma),\\
&   &  (-1)^{g(v)+\val(v)-1}(\tau u)^{r_1(v)}
       \prod_{(v,e)\in F(\Gamma)}d(e),\ \ \
       v\in V(\Gamma)^{(1)}.
\end{eqnarray*}
Recall that $r_1(v)=2g(v)-2+l(\mu(v))+\val(v)$ for $v\in V(\Gamma)^{(1)}$, and
$$
d_\Gamma^{(1)}=\left(\sum_{v\in V(\Gamma)^{(1)} }r_1(v)\right) -1
$$
is the virtual dimension of $\Mbar_\Gamma^{(1)}$.

We have
\begin{eqnarray*}
I^\Gamma(\tau)
&=&\frac{1}{|A_\Gamma|}\prod_{v\in V(\Gamma)^{(0)} }
   \left( |\mathrm{Aut}(\nu(v))|I^0_{g(v),\nu(v)}(p) \right)
   \left(\prod_{e\in E(\Gamma) } d(e)\right)^2\\
 &&  \cdot \prod_{ v \in V(\Gamma)^{(1)}}(-1)^{g(v)+\val(v)-1}
   \int_{[\Mbar_\Gamma^{(1)}]^{\mathrm{vir}} }
  \frac{(\tau u)^{d_\Gamma^{(1)}+1} }{-u-\psi^t}.
\end{eqnarray*}
Here $\nu(v)$ is the partition of $d(v)$ determined by
$\{ d(e):(v,e)\in F(\Gamma)\}$, as in Section \ref{graph}.

\smallskip

Recall that
$$
|A_\Gamma|=|\mathrm{Aut}(\Gamma)|
\prod_{e\in E}d(e),
$$
so
\begin{eqnarray*}
I^\Gamma(\tau)
& = & \frac{1}{|\mathrm{Aut}(\Gamma)|}
\prod_{v\in V(\Gamma)^{(0)} }
 (|\mathrm{Aut}(\nu(v))|I^0_{g(v),\nu(v)}(\tau))\\
&&\cdot \prod_{ v \in V(\Gamma)^{(1)}}\left( (-1)^{g(v)+\val(v)-1}
\prod_{(v,e)\in F(\Gamma)}d(e)\right)
\int_{[\Mbar_\Gamma^{(1)}]^{\mathrm{vir}} }
\frac{(\tau u)^{d_\Gamma^{(1)}+1}  }{-u-\psi^t} \\
& = & (-\tau)^{d_\Gamma^{(1)}+1}
\frac{1}{|\mathrm{Aut}(\Gamma)|}
\prod_{v\in V(\Gamma)^{(0)}}
(|\mathrm{Aut}(\nu(v))| I^0_{g(v),\nu(v)}(p) ) \\
& &\cdot \prod_{v\in V(\Gamma)^{(1)} }\left((-1)^{g(v)+\val(v)-1}
\prod_{(v,e)\in F}d(e)\right)
\int_{[\Mbar_\Gamma^{(1)}]^{\mathrm{vir}} } (\psi^t)^{d_\Gamma^{(1)}}\\
&=&\tau^{d_\Gamma^{(1)}+1}J^\Gamma(\tau).
\end{eqnarray*}

Note that $J^\Gamma(\tau)$ is a degree $r-d_\Gamma^{(1)}-1$ polynomial
in $\tau$, and
$$
J^\Gamma(-\tau-1)=(-1)^{d-l(\mu)+d_\Gamma^{(1)}+1}J^\Gamma(\tau).
$$

\subsection{Sum over graphs}\label{J1}
We have
$$
I^k_{g,\mu}(\tau)=\tau^k J^k_{g,\mu}(\tau),
$$
where $J^0_{g,\mu}(\tau)=I^0_{g,\mu}(\tau)$, and
$$
J^k_{g,\mu}(\tau)=\sum_{\Gamma\in G_{g,0}(\bP^1,\mu), \Gamma\neq\Gamma_0
d_\Gamma^{(1)}+1=k}J^\Gamma(\tau)
$$
for $k=1,\ldots,r$.
Note that $J^k_{g,\mu}(\tau)$ is a degree $r-k$ polynomial in $\tau$, and
$$
J^k_{g,\mu}(-\tau-1)=(-1)^{|\mu|-l(\mu)+k}J^k_{g, \mu}(\tau).
$$

For $k=1$, we have
$$
d_{\Gamma}^{(1)}=\sum_{v\in V(\Gamma)^{(1)} }r_1(v) -1 =0,
$$
so $V^S(\Gamma)^{(1)}$ consists of one single vertex $\bar{v}$
with $r_1(\bar{v})=1$. In particular, $g(\bar{v})=0$, and we have the
following two cases:
\begin{enumerate}
\item[Case 1:] $\mu(\bar{v})=(\mu_i)$, $\nu(\bar{v})=(j,k)$, where
 $j+k=\mu_i$. In this case, we say $\nu\in C(\mu)$ (cut) and define
 $a_{\mu,\nu}=jk$. Note that
\begin{eqnarray*}
&&\prod_{v\in V(\Gamma)^{(1)}}
\left((-1)^{g(v)+\val(v)-1}\prod_{(v,e)\in F}d(e)\right)
\int_{[\Mbar_\Gamma^{(1)}]^{\mathrm{vir}} } 1\\
&=&\left(-\mu_1\cdots\mu_{l(\mu)}\frac{jk}{\mu_i}\right)
    \frac{\mu_i}{\mu_1\cdots\mu_{l(\mu)}}\\
&=&-a_{\mu,\nu}.
\end{eqnarray*}
\item[Case 2:] $\mu(\bar{v})=(\mu_i,\mu_j)$, $\nu(\bar{v})=(\mu_i+\mu_j)$.
 In this case, we say $\nu\in J(\mu)$ (join) and define
 $b_{\mu,\nu}=\mu_i+\mu_j$. Note that
\begin{eqnarray*}
&&\prod_{v\in V(\Gamma)^{(1)} }
\left((-1)^{g(v)+\val(v)-1}\prod_{(v,e)\in F}d(e)\right)
\int_{[\Mbar_\Gamma^{(1)}]^{\mathrm{vir}} } 1\\
&=&\left(\mu_1\cdots\mu_{l(\mu)}\frac{\mu_i+\mu_j}{\mu_i\mu_j}\right)
    \frac{\mu_i\mu_j}{\mu_1\cdots\mu_{l(\mu)}}\\
&=&b_{\mu,\nu}.
\end{eqnarray*}

\end{enumerate}
So
\begin{eqnarray*}
J^1_{g,\mu}(\tau)
& = &
-\sum_{\nu\in J(\mu)}\frac{|\Aut(\nu)|}{|\mathrm{Aut}(\Gamma)|}
 b_{\mu,\nu}I^0_{g,\nu}(\tau)
+\sum_{\nu\in C(\mu)} \frac{|\Aut(\nu)|}{|\mathrm{Aut}(\Gamma)|}
 a_{\mu,\nu}I^0_{g-1,\nu}(\tau) \\
&& +\sum_{g_1+g_2=g,\nu^1\cup\nu^2\in C(\mu)}
 \frac{|\Aut(\nu^1)| |\Aut(\nu^2)|}{|\Aut(\Gamma)|}a_{\mu,\nu^1\cup\nu^2}
  I^0_{g_1,\nu^1}I^0_{g_2,\nu_2}(\tau) \\
& = &
-\sum_{\nu\in J(\mu)}I_1(\nu)I^0_{g,\nu}(\tau)
+\sum_{\nu\in C(\mu)} I_2(\nu)I^0_{g-1,\nu} (\tau)\\
&& +\sum_{g_1+g_2=g, \nu^1\cup\nu^2\in C(\mu)}
  I_3(\nu^1,\nu^2)
  I^0_{g_1,\nu^1}I^0_{g_2,\nu_2}(\tau),\\
\end{eqnarray*}
where $I_1$, $I_2$, and $I_3$ are defined as in \cite{Li-Zha-Zhe}.

We have
\begin{eqnarray*}
J_{g,\mu}^0(\tau)&=&\sqrt{-1}^{|\mu|-l(\mu)}\cC_{g,\mu}(\tau),\\
J_{g,\mu}^1(\tau)&=&\sqrt{-1}^{|\mu|-l(\mu)-1}\left(
\sum_{\nu\in J(\mu)}I_1(\nu) \cC_{g,\nu}(\tau)
+\sum_{\nu\in C(\mu)}I_2(\nu)\cC_{g-1,\nu}(\tau)\right.\\
& & \left. +\sum_{g_1+g_2=g,\nu^1\cup \nu^2\in C(\mu)}
   I_3(\nu^1,\nu^2) \cC_{g_1,\nu^1}(\tau) \cC_{g_2,\nu^2}(\tau)\right).
\end{eqnarray*}

It follows from the definition
that (\ref{eqn:CutJoin}) in Theorem~\ref{CCutJoin} is equivalent to
$$
\frac{d}{d\tau}J^0_{g,\mu}(\tau)= - J^1_{g,\mu}(\tau).
$$

\subsection{Final Calculations}
We have
\begin{eqnarray*}
F(\tau,x)
&=&\sum_{k=0}^r I_{g,\mu}^{r-k}(\tau)
x(x-1)\cdots(x-k+1)(x-k-1)\cdots (x-r)\\
&=&\sum_{k=0}^r \tau^{r-k} J^{r-k}_{g,\mu}(\tau)
x(x-1)\cdots(x-k+1)(x-k-1)\cdots (x-r)\\
&=&\sum_{k=0}^r \tau^k J^k_{g,\mu}(\tau)
x(x-1)\cdots(x-(r-k-1))(x-(r-k+1))\cdots (x-r).
\end{eqnarray*}
Therefore,
\begin{eqnarray*}
&& \mathrm{Br}_* e_T(V)\\
&=&\sum_{k=0}^r \tau^k J^k_{g,\mu}(\tau)
H(H-u)\cdots(H-(r-k-1)u)(H-(r-k+1)u)\cdots (H-ru).
\end{eqnarray*}

For $i=0,\ldots,r-1$, we have
\begin{eqnarray*}
&  & H^iH(H-u)\cdots(H-(r-k-1)u)(H-(r-k+1)u)\cdots (H-ru)\\
& = &((H-(r-k)u)+(r-k)u)^i H(H-u)\cdots(H-(r-k-1)u)\\
 & & \cdot(H-(r-k+1)u)\cdots (H-ru)\\
& = &((r-k)u)^i H(H-u)\cdots(H-(r-k-1)u)(H-(r-k+1)u)\cdots (H-ru)
\end{eqnarray*}
since
$$
H(H-u)\cdots(H-ru)=0.
$$
Therefore,
$$
\int_{\bP^r}\mathrm{Br}_* e_T(V)H^i
=u^i\sum_{k=0}^r(r-k)^i \tau^k J^k_{g,\mu}(\tau).
$$
Let $J^k_{g,\mu}(\tau)=\sum_{j=0}^{r-k}a_j^k \tau^j$. We have
$$
u^{-i}\int_{\bP^r}\mathrm{Br}_* e_T(V)H^i =
\sum_{l=0}^r\left(\sum_{j+k=l} (r-k)^i a_{j}^k\right)
\tau^l.
$$
Here is a crucial observation: as a polynomial in $\tau$,
$u^{-i}\int_{\bP^r}\mathrm{Br}_* e_T(V)H^i$
is of degree no more than $i$. Therefore,
$$
\sum_{j+k=l} (r-k)^i a_{j}^k=0
$$
for $0\leq i < l\leq r$. Now fix $l$ such that $1\leq l\leq r$.
We have
\begin{equation}\label{linear}
\sum_{k=0}^l (r-k)^i a_{l-k}^k =0,
\ \ \ \ 0\leq i<l,
\end{equation}
which is a system of $l$ linear equations of the
$l+1$ variables $\{a_{l-k}^k: k=0,\ldots,l\}$.

Both
$$
\{
(r-t)^i : i=0,\ldots,l-1
\}
$$
and
$$
\{
1, t, t(t-1), \ldots, t(t-1)\ldots(t-l+2)
\}
$$
are bases of the vector space
$$
\{ f(t)\in\bQ[t] :\deg(f)\leq l-1 \},
$$
so there exists an invertible $l\times l$ matrix
$(A_{ij})_{0\leq i,j\leq l-1}$ such that
$$
t(t-1)\cdots (t-i+1)=\sum_{j=0}^{l-1} A_{ij}(r-t)^j.
$$
In particular,
$$
k(k-1)\cdots (k-i+1)=\sum_{j=0}^{l-1} A_{ij}(r-k)^j.
$$
for $k=0,1,\ldots,l$, so (\ref{linear}) is equivalent to
$$
\sum_{k=0}^l k(k-1)\cdots (k-i+1) a_{l-k}^k =0,
\ \ \ \ 0\leq i<l,
$$
i.e.,
$$
\sum_{k=i}^l\frac{k!}{(k-i)!}a_{l-k}^k =0,
\ \ \ \ 0\leq i<l.
$$
The above equations can be rewritten as
$$
\left(
\begin{array}{ccccccc}
1 & 1      & \cdots & \cdots   & \cdots & \cdots & 1    \\
0 & 1!     & 2      & \cdots   & \cdots & \cdots & l   \\
0 & 0      & 2!     & 3\cdot 2 & \cdots & \cdots & l(l-1) \\
0 & 0      & 0      & 3!       & \cdots & \cdots & l(l-1)(l-2)\\
0 & \vdots & \vdots &\vdots    & \ddots & \vdots & \vdots\\
0 & \cdots & \cdots & \cdots   &0      & (l-1)! &l(l-1)\cdots 2
\end{array}\right)
\left(\begin{array}{c}a^0_l\\ a^1_{l-1}\\ \vdots\\ \vdots\\
   a^l_0 \end{array}\right)
= \left(\begin{array}{c}0\\\vdots\\ \vdots\\ \vdots\\ 0 \end{array}\right).
$$
The kernel is clearly one dimensional. One can check that the kernel
is given by
\begin{equation}\label{solution}
a_{l-k}^k=(-1)^k \frac{l!}{k!(l-k)!}a_l^0.
\end{equation}
Note that (\ref{solution}) for $l=1,\ldots,r$ is equivalent to
\begin{eqnarray} \label{eqn:Higher}
J_{g,\mu}^k(\tau)=\frac{(-1)^k}{k!}\frac{d^k}{d\tau^k}J_{g,\mu}^0(\tau)
\end{eqnarray}
for $k=0,\ldots,r$. In particular,
$$
J_{g,\mu}^1(\tau)=-\frac{d}{d\tau}J_{g,\mu}^0(\tau)
$$
which is equivalent to the cut-and-join equation (\ref{eqn:CutJoin})
in Theorem \ref{CCutJoin}.
Equation (\ref{eqn:Higher}) and the cut-and-join equation
imply that $J_{g,\mu}^k(\tau)$ can be obtained from $J_{g, \mu}^0(\tau)$
by repeating the cut-and-join operation $k$-times.

\section{Hurwitz Numbers} \label{hurwitz}

In this section, inspired by T. Graber and R. Vakil \cite{Gra-Vak1},
we use virtual localization on moduli spaces of relative stable
morphisms to recover some results on Hurwitz numbers.
We give a unified proof of the ELSV formula and the cut-and-join equation
for Hurwitz numbers.

\medskip

The Hirwitz numbers can be defined as 
\begin{equation}\label{eqn:hur}
H_{g,\mu}=\int_{[\Mbar_{g, 0}(\bP^1, \mu)]^{\mathrm{vir}}}
\mathrm{Br}^*{H^r}.
\end{equation}
We lift $H^r \in H^{2r}(\bP^r;\bZ)$ to
$\prod_{k=1}^r(H - w_k u)\in H^{2r}_T(\bP^r;\bZ)$, where $w_k\in \bZ$, and
compute
$$
H_{g,\mu}=\int_{[\Mbar_{g, 0}(\bP^1, \mu)]^{\mathrm{vir} }}
\mathrm{Br}^*\left(\prod_{k=1}^r(H - w_k u)\right)
$$
by virtual localization.

Let $p_k\in \bP^r$ be the $\bC^*$ fixed point defined as in Section \ref{branch},
and $f_k:p_k\to\bP^r$ be the inclusion. Then
$$
f_k^*(\prod_{l=1}^r(H - w_k u))=\left(\prod_{l=1}^r (k-w_l)\right)u^r\in H^{2r}_T(p_k).
$$
\begin{remark}
In the definition of $\Mbar_{g,0}(\bP^1,\mu)$ given in Section \ref{sec:MP},
if we order the $l(\mu)$ marked points on the domain as in \cite{Li1, Li2}, 
there will be an extra factor $1/|\Aut(\mu)|$ on the right-hand side of (\ref{eqn:hur}).
\end{remark}

\subsection{Contribution from each graph}
\subsubsection{The target is $\bP^1$}
Consider the graph $\Gamma^0\in \Mbar_{g,0}(\bP^1,\mu)$.
We have $\mathrm{Br}(F_{\Gamma^0})=\{p_r\}$.
We first consider the stable case.
By the Feynman rules derived in Appendix \ref{localization},
the contribution from $\Gamma^0$ is given by
$$\prod_{l=1}^r(r-w_l)\tilde{I}^0_{g,\mu},$$
where
\begin{eqnarray*}
\tilde{I}^0_{g,\mu}&=&
\frac{1}{|A_\Gamma|}\int_{\Mbar_{g,l(\mu)} }
\frac{u^r}{e(N^{\mathrm{vir}}_\Gamma)} \\
&=&\frac{1}{|\mathrm{Aut}(\mu)|}
   \prod_{i=1}^{l(\mu)}\frac{\mu_i^{\mu_i} }{\mu_i!}
    \int_{\Mbar_{g,l(\mu)} }
     \frac{\Lambda^\vee_{g}(u)u^{2g+2l(\mu)-3} }{
      \prod_{i=1}^{l(\mu)}(u-\mu_i\psi_i)}\\
&=&\frac{1}{|\mathrm{Aut}(\mu)|}
   \prod_{i=1}^{l(\mu)}\frac{\mu_i^{\mu_i} }{\mu_i!}
    \int_{\Mbar_{g,l(\mu)} }
     \frac{\Lambda^\vee_{g}(1)}{\prod_{i=1}^{l(\mu)}(1-\mu_i\psi_i)}
\end{eqnarray*}
In particular,
\begin{eqnarray*}
\tilde{I}^0_{0,\mu}
&=&\frac{1}{|\mathrm{Aut}(\mu)|}
     \prod_{i=1}^{l(\mu)}\frac{\mu_i^{\mu_i} }{\mu_i!}
    \int_{\Mbar_{0,l(\mu)} }
     \frac{1}{\prod_{i=1}^{l(\mu)}(1-\mu_i\psi_i)}\\
&=& \frac{1}{|\mathrm{Aut}(\mu)|}
  \left(\prod_{i=1}^{l(\mu)}\frac{\mu_i^{\mu_i} }{\mu_i!}
  \right) d^{l(\mu)-3}.
\end{eqnarray*}

The contribution from $\Gamma^{(0)}$ for the two unstable cases is also given by the
above formula:
\begin{eqnarray*}
      \tilde{I}^0_{0,(d)}
& = & \frac{d^d}{d!}\frac{1}{d^2}=\frac{d^{d-2}}{d!}.\\
       \tilde{I}^0_{0,(\mu_1,\mu_2)}
& = & \frac{1}{|\mathrm{Aut}((\mu_1,\mu_2))|}
       \frac{\mu_1^{\mu_1}\mu_2^{\mu_2}}{\mu_1!\mu_2!}\frac{1}{d}.
\end{eqnarray*}

\subsubsection{The target is $\bP^1[m]$, $m > 0$}
Consider $\Gamma\in G_{g,0}(\bP^1,\mu)$, $\Gamma\neq \Gamma^0$.
We have $\mathrm{Br}(F_\Gamma)=\{p_{r-d_\Gamma^{(1)}-1}\}$.
The contribution from $\Gamma$ is given by
$$
\prod_{l=1}^r(r-d_\Gamma^{(1)}-1-w_l)\tilde{I}^\Gamma,
$$
where
$$
\tilde{I}^\Gamma=
\frac{1}{|A_\Gamma|}\int_{[\Mbar_\Gamma]^{\mathrm{vir}} }
\frac{u^r}{e_T(N^{\mathrm{vir} }_\Gamma)}
=\frac{1}{|A_\Gamma|}\int_{[\Mbar_\Gamma]^{\mathrm{vir}} }
    \frac{\prod_{v\in V(\Gamma)} \tilde{B}_v }{-u-\psi^t},
$$
$$\tilde{B}_v=
\begin{cases}
 u^{2g(v)-2+\val(v)}, A_v\prod_{(v,e)\in F(\Gamma)}(A_e u^{d(e)})
& v\in V(\Gamma)^{(0)},\\
u^{r_1(v)} A_v,  &v\in V(\Gamma)^{(1)}.
\end{cases}
$$
More explicitly, in the notation of Appendix \ref{localization},
\begin{eqnarray*}
\tilde{B}_v
& = & \left(\prod_{(v,e)\in F(\Gamma)}d(e)\right)
\left(\prod_{(v,e)\in F(\Gamma)}\frac{d(e)^{d(e)}}{d(e)!}\right)
\frac{\Lambda^\vee_{g(v)}(u) u^{2g(v)+\val(v)-3} }{\prod_{(v,e)\in F(\Gamma)}
      \left(u-d(e)\psi_{(v,e)}\right)},
      \ \ \ v\in V^S(\Gamma)^{(0)},\\
&   & d(e)\frac{ d(e)^{d(e)-2} }{d(e)!},
      \ \ \ v\in V^I(\Gamma)^{(0)}, (v,e)\in F(\Gamma),\\
&   & d(e_1)d(e_2)\frac{d(e_1)^{d(e_1)}d(e_2)^{d(e)_2} }{d(e_1)!d(e_2)!}
      \frac{1}{d(v)},\ \ \
 v\in V^{II}(\Gamma)^{(0)}, (v,e_1),(v,e_2)\in F(\Gamma),\\
&   &  u^{r_1(v)}\prod_{(v,e)\in F(\Gamma)}d(e) ,\ \ \
        v\in V(\Gamma)^{(1)}.
\end{eqnarray*}

So
$$
\tilde{I}^\Gamma
= \frac{1}{|A_\Gamma|}
   \prod_{v\in V(\Gamma)^{(0)} }
   \left( |\mathrm{Aut}(\nu(v))|\tilde{I}^0_{g(v),\nu(v)} \right)
   \left(\prod_{e\in E(\Gamma) } d(e)\right)^2
   \int_{[\Mbar_\Gamma^{(1)}]^{\mathrm{vir}} }
\frac{u^{d_\Gamma^{(1)}+1} }{-u-\psi^t}.
$$

Recall that
$$
|A_\Gamma|=|\mathrm{Aut}(\Gamma)|
\prod_{e\in E}d(e)\prod_{v\in V^{II}(\Gamma)^{(1)} }d(v),
$$
so
$$
\tilde{I}^\Gamma=
\frac{(-1)^{d_\Gamma^{(1)}+1}}{|\mathrm{Aut}(\Gamma)|}
\prod_{v\in V(\Gamma)^{(0)}}(|\mathrm{Aut}(\nu(v))|I^0_{g(v),\nu(v)})
\left(\prod_{v\in V^S(\Gamma)^{(1)}}\prod_{(v,e)\in F}d(e)\right)
\int_{[\Mbar_\Gamma^{(1)}]^{\mathrm{vir}} }(\psi^t)^{d_\Gamma^{(1)}}.
$$

\subsection{Sum over Graphs}

For $k=1,\ldots,r$ define
$$
\tilde{I}^k_{g,\mu}=
\sum_{\Gamma\in G_{g,0}(\bP^1,\mu),\Gamma\neq \Gamma^0, d_\Gamma^{(1)}+1=k}
\tilde{I}^\Gamma.
$$
Then
$$
H_{g,\mu}=\sum_{k=0}^r \prod_{l=1}^r(r-k-w_l) \tilde{I}^k_{g,\mu}.
$$
for any $w_1,\ldots,w_r\in \bZ$.

Setting $(w_1,\ldots,w_r) =(0,1,\ldots,r-k-1,r-k+ 1,\ldots,r)$ yields
$$
H_{g,\mu}=(-1)^k k!(r-k)!\tilde{I}^k_{g,\mu}.
$$

The $k=0$ case implies the following formula due to
Ekedahl, Lando, Shapiro, and Vainshtein \cite{ELSV},
also proved by T. Graber and R. Vakil \cite{Gra-Vak1}:

\begin{equation}\label{eqn:ELSV}
H_{g,\mu} = r!\tilde{I}^0_{g,\mu}
=  \frac{r!}{|\mathrm{Aut}(\mu)|}
      \prod_{i=1}^{l(\mu)}\frac{ \mu_i^{\mu^i} }{\mu_i!}
      \int_{\Mbar_{g,l(\mu)} }
      \frac{\Lambda_g^\vee(1)}{\prod_{i=1}^{l(\mu)}(1-\mu_i\psi_i)}.
\end{equation}

For $k=1$, we have
$$
H_{g,\mu}=-(r-1)! \tilde{I}^1_{g,\mu},
$$
which is, by a derivation similar to that in Section \ref{J1}, equivalent to
\begin{equation}\label{eqn:k1}
\begin{split}
\frac{H_{g,\mu}}{(r-1)!}
=& \sum_{\nu\in J(\mu)}I_1(\nu)\tilde{I}^0_{g,\nu}
+\sum_{\nu\in C(\mu)} I_2(\nu)\tilde{I}^0_{g-1,\nu} \\
  & +\sum_{g_1+g_2=g}\sum_{\nu^1\cup\nu^2\in C(\mu)}
  I_3(\nu^1,\nu^2)\tilde{I}^0_{g_1,\nu^1}\tilde{I}^0_{g_2,\nu_2}.
\end{split}
\end{equation}

Combining (\ref{eqn:ELSV}) and (\ref{eqn:k1}), one obtains
\begin{eqnarray*}
&&\frac{H_{g,\mu}}{(r-1)!}=
\sum_{\nu\in J(\mu)}I_1(\nu)\frac{H_{g,\nu}}{(r-1)!}
+\sum_{\nu\in C(\mu)}I_2(\nu)\frac{H_{g-1,\nu} }{(r-1)!}\\
&& +\sum_{g_1+g_2=g}\sum_{\nu^1\cup\nu^2\in C(\mu)}
  I_3(\nu^1,\nu^2)
  \frac{H_{g_1,\nu^1} }{(2g_1-2+|\nu^1|+l(\nu^1))!}
  \frac{H_{g_2,\nu^2} }{(2g_2-2+|\nu^2|+l(\nu^2))!}
\end{eqnarray*}
which is equivalent to the cut-and-join equation
\begin{equation}
\begin{split}
&H_{g,\mu}
= \sum_{\nu\in J(\mu)}I_1(\nu)H_{g,\nu}+
\sum_{\nu\in C(\mu)}I_2(\nu)H_{g-1,\nu}\\
& +\sum_{g_1+g_2=g}\sum_{\nu^1\cup\nu^2\in C(\mu)}
\left(\begin{array}{cc} r-1 \\ 2g_1-2+|\nu^1|+l(\nu^1)\end{array} \right)
  I_3(\nu^1,\nu^2) H_{g_1,\nu^1} H_{g_2,\nu^2}.
\end{split}
\end{equation}

The above cut-and-join equation was first proved using combinatorics
by Goulden, Jackson and Vainstein in \cite{Gou-Jac-Vai} and later
proved using symplectic sum formula by
Li-Zhao-Zheng \cite{Li-Zha-Zhe} and Ionel-Parker
\cite[Section 15.2]{Ion-Par}.

\bigskip

{\bf Acknowledgements.}
We wish to thank Jun Li for explaining his work,
Jim Bryan, Bohui Chen, Tom Graber, Gang Liu, Ravi
Vakil for helpful conversations,
and Cumrun Vafa and Shing-Tung Yau for their interests in this work.
The first author wishes to thank the hospitality of IPAM where she
was a core participant of the Symplectic Geometry and Physics
Program and did most of her part of this work.
The second author is supported by an NSF grant.
The third author is partially supported by research
grants from NSFC and Tsinghua University.

\begin{appendix}
\section{ Localization} \label{localization}
In this appendix, we provide details of localization
in our particular case. Related results are discussed in \cite{Gra-Vak2}.
 We first introduce some notation.

Let $(w)$ be the 1-dimensional representation of $\bC^*$
given by $\lambda\cdot z =\lambda^w z$ for $\lambda\in\bC^*$,
$z\in\bC$. We do calculations on $\Mbar_\Gamma$
which is a finite cover of
$F_\Gamma\subset \Mbar_{g,0}(\bP^1,\mu)$,
so the weight $w$ of $\bC^*$-action can be fractional.

Given $\Gamma\in G_{g,0}(\bP^1,\mu)$, let
$$\left[f:(C,x_1,\ldots,x_{l(\mu)})\to\bP^1[m]\ \right]
$$
be a fixed point of the $\bC^*$-action on
$\Mbar_{g,0}(\bP^1,\mu)$ associated to
$\Gamma$. Given $(v,e)\in F(\Gamma)$, where
$v\in V^S(\Gamma)^{(0)}\cup V(\Gamma)^{(1)}$,
let $q_{(v,e)}\in C$ denote the node
at which $C_v$ and $C_e$ intersect. Let
$\psi_{(v,e)}$ denote the first Chern class
of the line bundle over $\Mbar_\Gamma$
whose fiber at
$$
\left[f:(C,x_1,\ldots,x_{l(\mu)})\to\bP^1[m]\ \right]
$$
is $T^*_{q_{(v,e)}}C_v$. Given $v\in V^{II}(\Gamma)^{(0)}$,
let $q_v$ denote the node at which $C_{e_1}$ and $C_{e_2}$
intersect, where $\{e_1,e_2\}=\{e\in E(\Gamma):(v,e)\in F(\Gamma)\}$.

\subsection{Virtual Normal Bundle}

The tangent space $T^1$ and the obstruction space $T^2$ of
$\Mbar_{g,0}(\bP^1,\mu)$
at
$$
\left[f: (C, x_1, \ldots, x_{l(\mu)}) \to \bP^1[m]\ \right]\in \Mbar_{g,0}(\bP^1,\mu)
$$
are given by the following two exact sequences \cite[Section 5.1]{Li2}:
\begin{eqnarray*}
0 &\to& \Ext^0(\Omega_C(D), \cO_C) \to H^0(\mathbf{D}^\bullet) \to T^1 \\
  &\to&  \Ext^1(\Omega_C(D), \cO_C) \to H^1(\mathbf{D}^\bullet) \to T^2 \to
  0\\
0 & \to &  H^0(C,f^*(\omega_{\bP^1[m]}(\log p_1^{(m)}))^\vee)
 \to H^0(\mathbf{D}^\bullet) \to
\oplus_{l=0}^{m-1} H_{\mathrm{et}}^0(\mathbf{R}_l^\bullet) \\
  & \to & H^1(C,f^*(\omega_{\bP^1[m]}(\log p_1^{(m)}))^\vee) \to H^1(\mathbf{D}^\bullet) \to
\oplus_{l=0}^{m-1} H_{\mathrm{et}}^1(\mathbf{R}_l^\bullet)\to 0
\end{eqnarray*}
where $D=x_1+\cdots+x_{l(\mu)}$, $\omega_{\bP^1[m]}$ is the dualizing sheaf of $\bP^1[m]$,
\begin{eqnarray*}
H^0_{\mathrm{et}}(\mathbf{R}_{l}^\bullet)&\cong&  \bigoplus_{q\in
f^{-1}(p_1^{(l)})}T_q \left(f^{-1}(\bP^1_{(l)})\right)\otimes
T^*_q \left(f^{-1}(\bP^1_{(l)})\right)\cong \bC^{\oplus n_l},\\
H^1_{\mathrm{et}}(\mathbf{R}_{l}^\bullet)&\cong&
(T_{p_1^{(l)}}\bP^1_{(l)}\otimes
T_{p_1^{(l)}}\bP^1_{(l+1)})^{\oplus(n_l-1)},
\end{eqnarray*}
and $n_l$ is the number of nodes over $p_1^{(l)}$.

Let $\tilde{f}=\pi[m]\circ f:C\to \bP^1$. Then
$$
f^*(\omega_{\bP^1[m]}(\log p_1^{(m)}))^\vee\cong \tilde{f}^*\cO_{\bP^1}(1).
$$

Let
$$
B_1 = \Ext^0(\Omega_C(D), \cO_C),\ \ \
B_2 = H^0(C,\tilde{f}^*\cO_{\bP^1}(1)),\ \ \
B_3 = \oplus_{l=0}^{m-1} H_{\mathrm{et}}^0(\mathbf{R}_l^\bullet),
$$
$$
B_4 = \Ext^1(\Omega_C(D), \cO_C),\ \ \
B_5 = H^1(C,\tilde{f}^*\cO_{\bP^1}(1)),\ \ \
B_6 = \oplus_{l=0}^{m-1} H_{\mathrm{et}}^1(\mathbf{R}_l^\bullet).
$$

We now assume that
$$
[f:(C,x_1,\ldots,x_{l(\mu)})\to \bP^1[m] ]\in F_\Gamma\subset\Mbar_{g,0}(\bP^1,\mu),
$$
where $\Gamma\in G_{g,0}(\bP^1,\mu)$, and $F_\Gamma$ is the set of fixed
points associated to $\Gamma$.
The $\bC^*$-action on $\Mbar_{g,0}(\bP^1,\mu)$ induces
$\bC^*$-actions on $B_1,\ldots,B_6$.
In particular, the $\bC^*$-actions on $B_2$ and $B_5$ come from
the $\bC^*$-action on $\cO_{\bP^1}(1)$ which acts on
$\cO_{\bP^1}(1)_{p_0}$ and $\cO_{\bP^1}(1)_{p_1}$
by $(1)$ and $(0)$, respectively.
Let $\tilde{B}_i$ denote the moving part of $B_i$ under the $\bC^*$-action. Then
each $\tilde{B}_i$ form a vector bundle over $\Mbar_\Gamma$. We will use
the same notation $\tilde{B}_i$ to denote the vector bundle.

Note that
$\tilde{B}_3=0$, and
$$
\tilde{B}_6=\left\{\begin{array}{ll}0,& m=0,\\
H^1_{\mathrm{et}}(\mathrm{R}_0^\bullet)=
(T_{p_1^{(0)}}\bP^1_{(0)}\otimes T_{p_1^{(0)}}\bP^1_{(1)})^{\oplus(n_0-1)}, &
m>0.
\end{array}\right.
$$
We have
$$
\frac{1}{e_T(N_\Gamma^{\mathrm{vir}} ) } = \frac{e_T(\tilde{T}^2)}{e_T(\tilde{T}^1)}
=\frac{e_T(\tilde{B}_1)e_T(\tilde{B}_5)e_T(\tilde{B}_6) }{e_T(\tilde{B}_2)e_T(\tilde{B}_4)}
$$
where $\tilde{T}^1$, $\tilde{T}^2$ are the moving parts of $T^1$, $T^2$, respectively.

\subsubsection{The target is $\bP^1$}
We have seen that there is only one graph $\Gamma^0$. Recall
that
$$
\Mbar_{\Gamma^0}=\left\{\begin{array}{ll}
\{\mathrm{point}\} & (g,l(\mu))= (0,1),(0,2),\\
\Mbar_{g,l(\mu)} &(g,l(\mu))\neq (0,1),(0,2).
\end{array}\right.
$$

We first compute $e_T(\tilde{B}_1)/e_T(\tilde{B}_4)$.
When $(g,l(\mu))\neq (0,1),(0,2)$.
The domain is
$$C=C_{v_0}\cup C_{e_1}\cup\cdots C_{e_{l(\mu)}}.$$
We have
$$
\tilde{B}_1=0,\ \ \ \tilde{B}_4=\bigoplus_{i=1}^{l(\mu)}
T_{q_i}C_{v_0}\otimes T_{q_i} C_{e_i},
$$
where $q_i=q_{(v_0,e_i)}$. So
$$
\frac{e_T(\tilde{B}_1)}{e_T(\tilde{B}_4)}=
\frac{1}{\prod_{i=1}^{l(\mu)}\left(\frac{u}{\mu_i}-\psi_i\right)},
$$
where  $\psi_i=\psi_{(v_0,e_i)}$.

When $(g,l(\mu))=(0,1)$, we have
$$
\tilde{B}_1=\left(\frac{1}{d}\right), \ \ \ \tilde{B}_4=0,
$$
so
$$
\frac{e_T(\tilde{B}_1)}{e_T(\tilde{B}_4)}=\frac{u}{d}.
$$

When  $(g,l(\mu))=(0,2)$, we have
$$
\tilde{B}_1=0, \ \ \ \tilde{B}_4=\left(\frac{1}{\mu_1}+\frac{1}{\mu_2}\right).
$$
so
$$
\frac{e_T(\tilde{B}_1)}{e_T(\tilde{B}_4)}=\frac{1}{\frac{u}{\mu_1}+\frac{u}{\mu_2}}.
$$

We next compute $e_T(\tilde{B}_5)/e_T(\tilde{B}_2)$. When
$(g,l(\mu))\neq (0,1), (0,2)$, consider the normalization sequence
$$
0 \to\tilde{f}^*\cO_{\bP^1}(1)\to
(\tilde{f}|_{C_{v_0}})^*\cO_{\bP^1}(1)\oplus
 \bigoplus_{i=1}^{l(\mu)} (\tilde{f}|_{C_{e_i}})^*\cO_{\bP^1}(1)\to
 \bigoplus_{i=1}^{l(\mu)}\cO_{\bP^1}(1)_{p_0}\to 0.
$$
The corresponding long exact sequence reads
\begin{eqnarray*}
0& \to & H^0(C,\tilde{f}^*\cO_{\bP^1}(1))\\
 & \to & H^0(C_{v_0},(\tilde{f}|_{C_{v_0}})^*\cO_{\bP^1}(1))\oplus
      \bigoplus_{i=1}^{l(\mu)} H^0(C_{e_i},(\tilde{f}|_{C_{e_i}})^*\cO_{\bP^1}(1))
   \to\bigoplus_{i=1}^{l(\mu)}\cO_{\bP^1}(1)_{p_0}\\
 &\to & H^1(C,\tilde{f}^*\cO_{\bP^1}(1))
  \to H^1(C_{v_0},(\tilde{f}|_{C_{v_0}})^*\cO_{\bP^1}(1))
   \oplus\bigoplus_{i=1}^{l(\mu)}
   H^1(C_{e_i},(\tilde{f}|_{C_{e_i}})^*\cO_{\bP^1}(1))\to 0.
\end{eqnarray*}
The representations of $\bC^*$ are
\begin{eqnarray*}
0& \to & H^0(C,\tilde{f}^*\cO_{\bP^1}(1))\to
 H^0(C_{v_0},\cO_{C_{v_0}})\otimes (1)\oplus\bigoplus_{i=1}^{l(\mu)}
\left(\bigoplus_{a=1}^{\mu_i}\left(\frac{a}{\mu_i}\right)\right)
\to \bigoplus_{i=1}^{l(\mu)}(1)\\
&\to&  H^1(C,\tilde{f}^*\cO_{\bP^1}(1))
\to  H^1(C_{v_0},\cO_{C_{v_0}})\otimes (1)\to 0
\end{eqnarray*}
So we have
\begin{equation}\label{Oone}
\frac{e_T(\tilde{B}_5)}{e_T(\tilde{B}_2)}
=\Lambda_g^\vee(u) u^{l(\mu)-1}
\prod_{i=1}^{l(\mu)}\left(\frac{\mu_i^{\mu_i}}{\mu_i !} u^{-\mu_i}\right).
\end{equation}
One can check that (\ref{Oone}) is also valid for
$(g,l(\mu))=(0,1),(0,2)$.

Finally, $\tilde{B}_6=0$, so $e_T(\tilde{B}_6)=1$.
We obtain the following Feynman rules:
$$
\frac{1}{e_T(N_{\Gamma^0}^{\mathrm{vir}}) }
=A_{v_0}\prod_{e\in E(\Gamma^0)} A_e,
$$
where
\begin{eqnarray*}
A_{v_0}
&=&\begin{cases} \frac{u}{d},&   (g,l(\mu))=(0,1)\\
\frac{u}{\frac{u}{\mu_1}+\frac{u}{\mu_2}}, & (g,l(\mu))=(0,2)\\
\frac{\Lambda_g^\vee(u)}{u}
\prod_{i=1}^{l(\mu)}\frac{u}{\frac{u}{\mu_i}-\psi_i}, &  (g,l(\mu))\neq (0,1),(0,2)
\end{cases}\\
A_e&=&\frac{d(e)^{d(e)}}{d(e)!} u^{-d(e)}.
\end{eqnarray*}

\subsubsection{The target is $\bP^1[m]$, $m>0$}
Let $\Gamma\in G_{g,0}(\bP^1,\mu)$, $\Gamma\neq \Gamma^0$.
We first compute
$e_T(\tilde{B}_1)/e_T(\tilde{B}_4)$.
We have
\begin{eqnarray*}
\tilde{B}_1 &=& \bigoplus_{v\in V^I(\Gamma)^{(0)}}\left(\frac{1}{d(v)}\right)\\
\tilde{B}_4 &=&\bigoplus_{v\in V^{II}(\Gamma)^{(0)} }
         \left(\sum_{(v,e)\in F(\Gamma)}\frac{1}{d(e)}\right)
         \oplus\bigoplus_{v\in V^S(\Gamma)^{(0)} }
           \left(\bigoplus_{(v,e)\in F(\Gamma)}
          T_{q_{(v,e)}}C_v\otimes T_{q_{(v,e)}}C_e\right)\\
      &&  \oplus\bigoplus_{v\in V(\Gamma)^{(1)} }
           \left(\bigoplus_{(v,e)\in F(\Gamma)}
           T_{q_{(v,e)}}C_v\otimes T_{q_{(v,e)}}C_e\right)
\end{eqnarray*}
So
\begin{eqnarray*}
\frac{e_T(\tilde{B}_1)}{e_T(\tilde{B}_4)}
&=&\prod_{v\in V^I(\Gamma)^{(0)}}\frac{u}{d(v)}
\prod_{v\in V^{II}(\Gamma)^{(0)} }
   \left( \sum_{(v,e)\in F(\Gamma)}\frac{u}{d(e)} \right)^{-1}\\
&&\cdot\prod_{v\in V^S(\Gamma)^{(0)} }
           \left(\sum_{(v,e)\in F(\Gamma)}
          \frac{1}{\frac{u}{d(e)}-\psi_{(v,e)}}\right)
       \prod_{v\in V(\Gamma)^{(1)} }
        \left(\prod_{(v,e)\in F(\Gamma)}
          \frac{1}{\frac{-u}{d(e)}-\psi_{(v,e)}}\right)
\end{eqnarray*}

We next compute $e_T(\tilde{B}_5)/e_T(\tilde{B}_2)$. Consider the
normalization sequence
\begin{eqnarray*}
0 &\to&
\tilde{f}^*\cO_{\bP^1}(1)\to
\bigoplus_{v\in V^{S}(\Gamma)^{(0)}\cup V(\Gamma)^{(1)} }
(\tilde{f}|_{C_v})^*\cO_{\bP^1}(1)
\oplus  \bigoplus_{e\in E(\Gamma)}(\tilde{f}|_{C_e})^*\cO_{\bP^1}(1)\\
  &\to& \bigoplus_{v\in V^{II}(\Gamma)^{(0)}} \cO_{\bP^1}(1)_{p_0}
\oplus\bigoplus_{v\in V^S(\Gamma)^{(0)}}\left(\bigoplus_{(v,e)\in F}
     \cO_{\bP^1}(1)_{p_0}\right)
\oplus\bigoplus_{v\in V(\Gamma)^{(1)}}\left(\bigoplus_{(v,e)\in F}
     \cO_{\bP^1}(1)_{p_1}\right)\to 0.
\end{eqnarray*}
The corresponding long exact sequence reads
\begin{eqnarray*}
0& \to & H^0(C,\tilde{f}^*\cO_{\bP^1}(1))\\
 & \to & \bigoplus_{v\in V^{S}(\Gamma)^{(0)}\cup V(\Gamma)^{(1)} }
         H^0(C_v,(\tilde{f}|_{C_v})^*\cO_{\bP^1}(1))
         \oplus \bigoplus_{e\in E(\Gamma)}
 H^0(C_e,(\tilde{f}|_{C_e})^*\cO_{\bP^1}(1))\\
&\to&\bigoplus_{v\in V^{II}(\Gamma)^{(0)}} \cO_{\bP^1}(1)_{p_0}
\oplus\bigoplus_{v\in V^S(\Gamma)^{(0)}}\left(\bigoplus_{(v,e)\in F}
     \cO_{\bP^1}(1)_{p_0}\right)
\oplus\bigoplus_{v\in V(\Gamma)^{(1)}}\left(\bigoplus_{(v,e)\in F}
     \cO_{\bP^1}(1)_{p_1}\right) \\
 &\to & H^1(C,\tilde{f}^*\cO_{\bP^1}(1))\\
 & \to & \bigoplus_{v\in V^{S}(\Gamma)^{(0)}\cup V(\Gamma)^{(1)} }
           H^0(C_v,(\tilde{f}|_{C_v})^*\cO_{\bP^1}(1))
          \oplus \bigoplus_{e\in E(\Gamma)}
 H^0(C_e,(\tilde{f}|_{C_e})^*\cO_{\bP^1}(1))\to 0.
\end{eqnarray*}
The representations of $\bC^*$ are
\begin{eqnarray*}
0& \to & H^0(C,\tilde{f}^*\cO_{\bP^1}(1))\\
 & \to & \bigoplus_{v\in V^{S}(\Gamma)^{(0)} }
         H^0(C_v,\cO_{C_v})\otimes (1)
  \oplus \bigoplus_{v\in V(\Gamma)^{(1)} }
         H^0(C_v,\cO_{C_v})\otimes (0)
  \oplus \bigoplus_{e\in E(\Gamma)}
\left(\bigoplus_{a=1}^{d(e)}\left(\frac{a}{d(e)}\right)\right) \\
&\to&\bigoplus_{v\in V^{II}(\Gamma)^{(0)} } (1)
\oplus\bigoplus_{v\in V^S(\Gamma)^{(0)} }\left(\bigoplus_{(v,e)\in F}
     (1)\right)
\oplus\bigoplus_{v\in V(\Gamma)^{(1)} }\left(\bigoplus_{(v,e)\in F}
     (0)\right) \\
 &\to & H^1(C,\tilde{f}^*\cO_{\bP^1}(1))\\
 &\to & \bigoplus_{v\in V^{S}(\Gamma)^{(0)} }
         H^1(C_v,\cO_{C_v})\otimes (1)
  \oplus  \bigoplus_{v\in  V(\Gamma)^{(1)} }
H^1(C_v,\cO_{C_v})\otimes (0) \to 0
\end{eqnarray*}
So
$$
\frac{e_T(\tilde{B}_5)}{e_T(\tilde{B}_2)}=
\prod_{v\in V(\Gamma)^{(0)}}
\left(\Lambda_{g(v)}^\vee(u)u^{\val(v)-1}\right)
\prod_{e\in E(\Gamma)}\left(\frac{d(e)^{d(e)}}{d(e)!} u^{-d(e)}\right).
$$
Finally,
$$
\tilde{B}_6=\left(T_{p_1^{(0)}}\bP^1_{(0)}
             \otimes T_{p_1^{(0)}}\bP^1_{(1)}\right)^{|E(\Gamma)|-1},
$$
so
$$
e_T(\tilde{B}_6)=(-u-\psi^t)^{|E(\Gamma)|-1},
$$
where $\psi^t$ is the first Chern class of the line bundle
over $\Mbar_\Gamma^{(1)}$ whose fiber at
$$[\hat{f}:\hat{C}\to\bP^1(m)]$$
is $T^*_{p_1^{(0)}}\bP^1(m)$.
Note that $\psi^t=d(e)\psi_{(v,e)}$ for $v\in V(\Gamma)^{(1)}$,
$(v,e)\in F$.

\medskip
We have the following Feynman rules:
$$\frac{1}{e_T(N_{\Gamma}^{\mathrm{vir} })}
= \frac{1}{-u-\psi^t}\prod_{v \in V(\Gamma)} A_v
\prod_{e \in E(\Gamma)} A_e,
$$
where
\begin{eqnarray*}
A_v &= &\left\{
\begin{array}{ll}
\prod_{(v,e)\in F(\Gamma)}\frac{-u-\psi^t}{\frac{-u}{d(e)}-\psi_{(v,e)} }
 =\prod_{(v,e)\in F(\Gamma)} d(e), &
 v \in V(\Gamma)^{(1)}\\
\frac{\Lambda_{g(v)}^{\vee}(u)}{u}\prod_{(v,e)\in F(\Gamma)}
\frac{u}{\frac{u}{d(e)} - \psi_{(v, e)}},
   & v \in V^S(\Gamma)^{(0)}, \\
\frac{u}{d(v)}, & v \in V^I(\Gamma)^{(0)}, \\
\frac{1}{\frac{1}{d(e_1)} + \frac{1}{d(e_2)}},
   & v \in V^{II}(\Gamma)^{(0)},
(v, e_1), (v, e_2) \in F(\Gamma),
\end{array} \right.\\
 A_e&=& \frac{d(e)^{d(e)} }{d(e)!}u^{-d(e)}.
\end{eqnarray*}

The degree of $e_T(N^{\mathrm{vir}}_{\Gamma})$ is
\begin{eqnarray*}
 & & 1+ \sum_{v \in V^S(\Gamma)^{(0)} } (1 - g(v)) -|V^I(\Gamma)^{(0)}|
     + \sum_{e \in E(\Gamma)} d(e)\\
& = & 1+ \sum_{v \in V^S(\Gamma)^{(0)} } (1 - g(v)) -|V^I(\Gamma)^{(0)}|+ d.
\end{eqnarray*}
We have seen that the virtual dimension of $F_\Gamma$ is
$$
2g-3+l(\mu)+\sum_{v \in V^S(\Gamma)^{(0)} }(g(v)-1) + |V^I(\Gamma)^0|.
$$
We have
\begin{eqnarray*}
& &2g-3+l(\mu)+\sum_{v \in V^S(\Gamma)^{(0)} }(g(v)-1) + |V^I(\Gamma)^0|  \\
& & + 1+ \sum_{v \in V^S(\Gamma)^{(0)} } (1 - g(v)) -|V^I(\Gamma)^{(0)}|+d \\
&=& 2g-2+d+l(\mu)
\end{eqnarray*}
as expected.

\subsection{The bundle $V_{D}$}\label{D}
The short exact sequence
$$
0\to\cO_C(-D)\to \cO_C \to \cO_D\to 0
$$
gives rise to a long exact sequence
\begin{eqnarray*}
0&\to& H^0(C,\cO_C(-D) )\to H^0(C,\cO_C)\to
\bigoplus_{i=1}^{l(\mu)}\cO_{x_i}\\
&\to& H^1(C,\cO_C(-D) )\to H^1(C,\cO_C)\to 0
\end{eqnarray*}
The representations of $\bC^*$ are
$$
0\to (\tau)\to \bigoplus_{i=1}^{l(\mu)}(\tau) \to
H^1(C,\cO_C(-D))\otimes (\tau)\to H^1(C,\cO_C)\otimes(\tau)\to 0.
$$
So
$$
e_T(V_D)=e_T(V_0)(\tau u)^{l(\mu)-1},
$$
where
$$V_0=R^1\pi_*\cO_{\mathcal{U}_{g,\mu}}.$$
Recall that $\pi:\mathcal{U}_{g,\mu}\to \Mbar_{g,0}(\bP^1,\mu)$ is
the universal curve.

\subsubsection{The target is $\bP^1$}
There is only one graph $\Gamma^0$ in this case.
When $(g,l(\mu))\neq (0,1), (0,2)$,
consider the normalization sequence
$$
0\to \cO_C\to\cO_{C_{v_0}}\oplus \bigoplus_{i=1}^{l(\mu)} \cO_{C_{e_i}}
 \to \bigoplus_{i=1}^{l(\mu)}\cO_{q_i}\to 0.
$$
The corresponding long exact sequence reads
\begin{eqnarray*}
0& \to &  H^0(C,\cO_C)
   \to H^0(C_{v_0},\cO_{C_{v_0}})\oplus
      \bigoplus_{i=1}^{l(\mu)} H^0(C_{e_i},\cO_{C_{e_i}})
   \to\bigoplus_{i=1}^{l(\mu)}\cO_{q_i}\\
 &\to & H^1(C,\cO_C)
  \to H^1(C_{v_0},\cO_{C_{v_0}})
   \oplus\bigoplus_{i=1}^{l(\mu)}
   H^1(C_{e_i},\cO_{C_{e_i}})\to 0
\end{eqnarray*}
The representations of $\bC^*$ are
\begin{eqnarray*}
0& \to&   (\tau)
   \to (\tau)\oplus
      \bigoplus_{i=1}^{l(\mu)} (\tau)
   \to\bigoplus_{i=1}^{l(\mu)} (\tau)\\
& \to & H^1(C,\cO_C)\otimes(\tau)
  \to H^1(C_{v_0},\cO_{C_{v_0}})\otimes(\tau)\to 0.
\end{eqnarray*}
So
$$
i_\Gamma^* e_T(V_0)=\Lambda_g^\vee(\tau u)
$$
which is clearly also valid for $(g,l(\mu))=(0,1),(0,2)$.

We have
$$
i_{\Gamma^0}^* e_T(V_D)=A^D_{v_0},
$$
where
$$
A_{v_0}^D = \Lambda_g^{\vee}(\tau u) \cdot
(\tau u)^{l(\mu)-1}.
$$

Note that the degree of $i^*_{\Gamma^0} e_T(V_D)$ is $l(\mu)+g -
1$, as expected.

\subsubsection{The target is $\bP^1[m]$, $m>0$}
Given a graph $\Gamma\in G_{g,0}(\bP^1,\mu)$, $\Gamma\neq \Gamma^0$,
consider the normalization sequence
\begin{eqnarray*}
0 &\to&\cO_C\to
\bigoplus_{v\in V^{S}(\Gamma)^{(0)}\cup V(\Gamma)^{(1)} }\cO_{C_v}
 \oplus \bigoplus_{e\in E(\Gamma)}\cO_{C_e}\\
&\to&\bigoplus_{v\in V^{II}(\Gamma)^{(0)}} \cO_{q_v}
\oplus\bigoplus_{v\in V^S(\Gamma)^{(0)}\cup V(\Gamma)^{(1)}}\left(\bigoplus_{(v,e)\in F}
     \cO_{q_{(v,e)}}\right)\to 0.
\end{eqnarray*}
The corresponding long exact sequence reads
\begin{eqnarray*}
0& \to & H^0(C,\cO_C)
\to  \bigoplus_{v\in V^{S}(\Gamma)^{(0)}\cup V(\Gamma)^{(1)} }
         H^0(C_v,\cO_{C_v})
         \oplus \bigoplus_{e\in E(\Gamma)}
 H^0(C_e,\cO_{C_e})\\
&\to&\bigoplus_{v\in V^{II}(\Gamma)^{(0)}} \cO_{q_v}
\oplus\bigoplus_{v\in V^S(\Gamma)^{(0)}\cup V(\Gamma)^{(1)}}\left(\bigoplus_{(v,e)\in F}
     \cO_{q_{(v,e)}}\right)\\
 &\to & H^1(C,\cO_C)
\to \bigoplus_{v\in V^{S}(\Gamma)^{(0)}\cup V(\Gamma)^{(1)} }
           H^1(C_v,\cO_{C_v})
          \oplus \bigoplus_{e\in E(\Gamma)}
 H^1(C_e,\cO_{C_e})\to 0
\end{eqnarray*}
The representations of $\bC^*$ are
\begin{eqnarray*}
0& \to &(\tau)
\to \bigoplus_{v\in V^{S}(\Gamma)^{(0)}\cup V(\Gamma)^{(1)} }(\tau)
     \oplus \bigoplus_{e\in E(\Gamma)}(\tau)\\
&\to&\bigoplus_{v\in V^{II}(\Gamma)^{(0)}} (\tau)
\oplus\bigoplus_{v\in V^S(\Gamma)^{(0)}\cup V(\Gamma)^{(1)}}\left(\bigoplus_{(v,e)\in F}
     (\tau)\right)\\
 &\to & H^1(C,\cO_C)\otimes(\tau)
\to \bigoplus_{v\in V^{S}(\Gamma)^{(0)}\cup V(\Gamma)^{(1)} }
           H^1(C_v,\cO_{C_v})\otimes(\tau)\to 0
\end{eqnarray*}
So
\begin{eqnarray*}
i_\Gamma^* e_T(V_0)&=&(\tau u)^{|E(\Gamma)|-|V(\Gamma)|+1}
\prod_{v\in V(\Gamma)}\Lambda_{g(v)}^\vee(\tau u)\\
i_\Gamma^* e_T(V_D)&=&(\tau u)^{|E(\Gamma)|-|V(\Gamma)|+l(\mu)}
\prod_{v\in V(\Gamma)}\Lambda_{g(v)}^\vee(\tau u)\\
\end{eqnarray*}

We have the following Feynman rules:
$$i_{\Gamma}^*e_T(V_D) = \prod_{v\in V(\Gamma)} A_v^D,$$
where
$$
 A_v^{D} = \begin{cases}
\Lambda_{g(v)}^{\vee}(\tau u) \cdot (\tau u)^{\val(v)-1},
& v \in V(\Gamma)^{(0)}, \\
\Lambda_{g(v)}^{\vee}(\tau u) \cdot
(\tau u)^{l(\mu(v))-1}
& v \in V(\Gamma)^{(1)},
\end{cases} \\
$$

Note that the degree of $i_{\Gamma}^*e_T(V_D)$ is
\begin{eqnarray*}
&   & \sum_{v\in V(\Gamma)^{(0)}}(g(v)+\val(v)-1)
      +\sum_{v\in V(\Gamma)^{(1)}}(g(v)+l(\mu(v))-1)\\
& = & \sum_{v\in V(\Gamma)}g(v)+ |E(\Gamma)|-|V(\Gamma)|+l(\mu(v))\\
& = & l(\mu) + \left(\sum_{v\in V(\Gamma)}g(v) + |E(\Gamma)|-|V(\Gamma)|+1\right)-1\\
& = & l(\mu) + g - 1
\end{eqnarray*}
as expected.

\subsection{The bundle $V_{D_d}$} \label{Dd}
\subsubsection{The target is $\bP^1$}
There is only one graph $\Gamma^0$ in this case.
When $(g,l(\mu))\neq (0,1), (0,2)$,
consider the normalization sequence
$$
0\to\tilde{f}^*\cO_{\bP^1}(-1)\to
(\tilde{f}|_{C_{v_0}})^*\cO_{\bP^1}(-1)\oplus
 \bigoplus_{i=1}^{l(\mu)}(\tilde{f}|_{C_{e_i}})^*\cO_{\bP^1}(-1)
\to\bigoplus_{i=1}^{l(\mu)}\cO_{\bP^1}(-1)_{p_0}\to 0.
$$
The corresponding long exact sequence reads
\begin{eqnarray*}
0& \to & H^0(C,\tilde{f}^*\cO_{\bP^1}(-1))\\
 & \to & H^0(C_{v_0},(\tilde{f}|_{C_{v_0}})^*\cO_{\bP^1}(-1))\oplus
      \bigoplus_{i=1}^{l(\mu)} H^0(C_{e_i},(\tilde{f}|_{C_{e_i}})^*\cO_{\bP^1}(-1))
   \to\bigoplus_{i=1}^{l(\mu)}\cO_{\bP^1}(-1)_{p_0}\\
 &\to & H^1(C,\tilde{f}^*\cO_{\bP^1}(-1))\\
 &\to & H^1(C_{v_0},(\tilde{f}|_{C_{v_0}})^*\cO_{\bP^1}(-1))
   \oplus\bigoplus_{i=1}^{l(\mu)}
   H^1(C_{e_i},(\tilde{f}|_{C_{e_i}})^*\cO_{\bP^1}(-1))\to 0.
\end{eqnarray*}
The representations of $\bC^*$ are
\begin{eqnarray*}
0& \to & 0
  \to  H^0(C_{v_0},\cO_{C_{v_0}})\otimes (-\tau-1)
  \to \bigoplus_{i=1}^{l(\mu)}(-\tau-1)\\
&\to& H^1(C,\tilde{f}^*\cO_{\bP^1}(-1))\\
&\to&
 H^1(C_{v_0},\cO_{C_{v_0}})\otimes (-\tau-1)
 \oplus\bigoplus_{i=1}^{l(\mu)}
\left(\bigoplus_{a=1}^{d(e)-1}\left(-\tau-\frac{a}{d(e)} \right)\right)\to 0
\end{eqnarray*}

We have the following Feynman rules:
$$i_{\Gamma^0}^*e_T(V_{D_d}) = A_{v_0}^{D_d}
\prod_{e\in E(\Gamma)} A_e^{D_d},$$
where
\begin{eqnarray*}
A_{v_0}^{D_d} &=& \Lambda_g^{\vee}((-\tau-1) u) \cdot
((-\tau-1)u)^{l(\mu)-1}\\
A_e^{D_d}&=&
\frac{\prod_{a=1}^{d(e)-1} (d(e)\tau + a)}{d(e)^{d(e)-1}}(-u)^{d(e)-1}.
\end{eqnarray*}
It is easily checked that the above Feynmann rules are also valid for
$(g,l(\mu))=(0,1),(0,2)$.

Note that the degree of $i_\Gamma^*e_T(V_{D_d})$ is
$d+g-1$, as expected.

\subsubsection{The target is $\bP^1[m]$, $m>0$}
Given a graph $\Gamma\in G_{g,0}(\bP^1,\mu)$, $\Gamma\neq \Gamma^0$,
consider the normalization sequence
\begin{eqnarray*}
0 & \to & \tilde{f}^*\cO_{\bP^1}(-1)\to
 \bigoplus_{v\in V^{S}(\Gamma)^{(0)}\cup V(\Gamma)^{(1)} }
(\tilde{f}|_{C_v})^*\cO_{\bP^1}(-1)
 \oplus \bigoplus_{e\in E(\Gamma)}(\tilde{f}|_{C_e})^*\cO_{\bP^1}(-1)\\
&\to&\bigoplus_{v\in V^{II}(\Gamma)^{(0)}} \cO_{\bP^1}(-1)_{p_0}
\oplus\bigoplus_{v\in V^S(\Gamma)^{(0)}}\left(\bigoplus_{(v,e)\in F}
     \cO_{\bP^1}(-1)_{p_0}\right)\\
&&\oplus\bigoplus_{v\in V(\Gamma)^{(1)}}\left(\bigoplus_{(v,e)\in F}
     \cO_{\bP^1}(-1)_{p_1}\right)\to 0.
\end{eqnarray*}
The corresponding long exact sequence reads
\begin{eqnarray*}
0& \to & H^0(C,\tilde{f}^*\cO_{\bP^1}(-1))\\
 & \to & \bigoplus_{v\in V^{S}(\Gamma)^{(0)}\cup V(\Gamma)^{(1)} }
         H^0(C_v,(\tilde{f}|_{C_v})^*\cO_{\bP^1}(1))
         \oplus \bigoplus_{e\in E(\Gamma)}
 H^0(C_e,(\tilde{f}|_{C_e})^*\cO_{\bP^1}(-1))\\
&\to&\bigoplus_{v\in V^{II}(\Gamma)^{(0)}} \cO_{\bP^1}(-1)_{p_0}
\oplus\bigoplus_{v\in V^S(\Gamma)^{(0)}}\left(\bigoplus_{(v,e)\in F}
     \cO_{\bP^1}(-1)_{p_0}\right)\\
&& \oplus\bigoplus_{v\in V(\Gamma)^{(1)}}\left(\bigoplus_{(v,e)\in F}
     \cO_{\bP^1}(-1)_{p_1}\right) \\
 &\to & H^1(C,\tilde{f}^*\cO_{\bP^1}(-1))\\
 & \to & \bigoplus_{v\in V^{S}(\Gamma)^{(0)}\cup V(\Gamma)^{(1)} }
           H^1(C_v,(\tilde{f}|_{C_v})^*\cO_{\bP^1}(1))
          \oplus \bigoplus_{e\in E(\Gamma)}
 H^1(C_e,(\tilde{f}|_{C_e})^*\cO_{\bP^1}(-1))\to 0.
\end{eqnarray*}
The representations of $\bC^*$ are
\begin{eqnarray*}
0& \to & H^0(C,\tilde{f}^*\cO_{\bP^1}(-1))\\
&\to& \bigoplus_{v\in V^{S}(\Gamma)^{(0)}}H^0(C_v,\cO_{C_v})\otimes (-\tau-1)
  \oplus \bigoplus_{v\in V(\Gamma)^{(1)}}H^0(C_v,\cO_{C_v})\otimes (-\tau) \\
 &\to&\bigoplus_{v\in V^{II}(\Gamma)^{(0)}} (-\tau-1)
\oplus \bigoplus_{v\in V^S(\Gamma)^{(0)}}\left(\bigoplus_{(v,e)\in F}
     (-\tau-1)\right)
\oplus \bigoplus_{v\in V(\Gamma)^{(1)}}\left(\bigoplus_{(v,e)\in F}
     (-\tau)\right) \\
 &\to & H^1(C,\tilde{f}^*\cO_{\bP^1}(-1))\\
&\to& \bigoplus_{v\in V^{S}(\Gamma)^{(0)}}H^1(C_v,\cO_{C_v})\otimes (-\tau-1)
  \oplus \bigoplus_{v\in V(\Gamma)^{(1)}}H^1(C_v,\cO_{C_v})\otimes (-\tau)\\
&&  \oplus \bigoplus_{e\in E(\Gamma)}
\left(\prod_{a=1}^{d(e)-1}\left(-\tau-\frac{a}{d(e)}\right) \right)\to 0
\end{eqnarray*}

We have the following Feynman rules:
$$i_{\Gamma}^*e_T(V_{D_d}) = \prod_v A_v^{D_d} \cdot
\prod_{e\in E(\Gamma)} A_e^{D_d},$$
where
\begin{eqnarray*}
&& A_v^{D_d} = \begin{cases}
\Lambda_{g(v)}^{\vee}((-\tau-1) u) \cdot
((-\tau-1)u)^{\val(v)-1}, & v \in V(\Gamma)^{(0)}, \\
\Lambda_{g(v)}^{\vee}(-\tau u) \cdot
(-\tau u)^{\val(v)-1}, & v \in V(\Gamma)^{(1)}, \\
\end{cases} \\
&& A_e^{D_d}
=\frac{\prod_{a=1}^{d(e)-1} (d(e) \tau + a)}{d(e)^{d(e)-1}}(-u)^{d(e)-1}.
\end{eqnarray*}

Note that the degree of $i_{\Gamma}^*e_T(V_{D_d})$ is
\begin{eqnarray*}
&   & \sum_{v\in V(\Gamma)}(g(v)+\val(v)-1)+\sum_{e\in E(\Gamma)}(d(e)-1)\\
& = & \sum_{v\in V(\Gamma)}g(v) + 2|E(\Gamma)|-|V(\Gamma)|+d-|E|\\
& = & d +  \left(\sum_{v\in V(\Gamma)}g(v)+|E(\Gamma)|-|V(\Gamma)|+1\right) -1\\
& = & d + g - 1
\end{eqnarray*}
as expected.

\subsection{The obstruction bundle}
Combining results in Section \ref{Dd} and Section \ref{D}, we
obtain Feynman rules for the obstruction bundle
$$V=V_D\oplus V_{D_d}.$$
\subsubsection{The target is $\bP^1$}
We have
$$i_{\Gamma^0}^*(e_T(V)) = A_{v_0}^V \prod_{e\in E(\Gamma^0)} A_e^V,$$
where
\begin{eqnarray*}
A_{v_0}^V = A_{v_0}^D A_{v_0}^{D_d}&=&
\Lambda_g^{\vee}(\tau u)\Lambda_g^{\vee}(-(\tau+1)u) \cdot
(\tau u(-\tau-1)u)^{l(\mu)-1}\\
A_e^V=A_e^D A_e^{D_d}&=&
\frac{\prod_{a=1}^{d(e)-1} (d(e) p + a)}{d(e)^{d(e)-1}}(-u)^{d(e)-1}.
\end{eqnarray*}

\subsubsection{The target is $\bP^1[m]$, $m>0$}
For any $\Gamma\in G_{g,0}(\bP^1,\mu)$, $\Gamma\neq \Gamma^0$, we have
$$i_{\Gamma}^*(e_T(V)) =
\prod_{v\in V(\Gamma)} A_v^V \cdot \prod_{e\in E(\Gamma)} A_e^V,$$
where
\begin{eqnarray*}
&& A_v^V = A_v^D A_v^{D_d}= \begin{cases}
\Lambda_{g(v)}^{\vee}(\tau u)\Lambda_{g(v)}^{\vee}(-(\tau+1)u) \cdot
(\tau u(-\tau-1)u)^{\val(v)-1}, & v \in V(\Gamma)^{(0)}, \\
(-1)^{g(v)+\val(v)-1}
(\tau u)^{r_1(v)}, &   v \in V(\Gamma)^{(1)}
\end{cases} \\
&& A_e^V=A_e^D A_e^{D_d}
= \frac{\prod_{a=1}^{d(e)-1} (d(e) \tau + a)}{d(e)^{d(e)-1}} (-u)^{d(e)-1}.
\end{eqnarray*}
Recall that $r_1(v)=2g(v)-2+l(\mu(v))+\val(v)$ for $v\in V(\Gamma)^{(1)}$.
\end{appendix}

\end{document}